\documentclass[a4paper,10pt]{article}
\usepackage[T1]{fontenc}
\usepackage{graphicx,amsmath,amsfonts,amssymb,amsthm,framed,psfrag,
color,cancel,dsfont,textcomp
}
\usepackage{mathrsfs,eufrak}
\usepackage{algorithm2e,tikz}
\usetikzlibrary{arrows.meta}
\tikzset{%
  >={Latex[width=2mm,length=2mm]},
            base/.style = {rectangle, rounded corners, draw=black,
                           minimum width=3.7cm, minimum height=1.2cm,
                           text centered, font=\sffamily},
  activityStarts/.style = {base, fill=blue!15},
       startstop/.style = {base, fill=red!15},
    activityRuns/.style = {base, fill=green!15},
         process/.style = {base, minimum width=5.3cm, minimum height=1.5cm, fill=orange!10,
                           font=\ttfamily},
}

\newcommand{\bre}{{\mbox{\sc ExitBm}}}
\newcommand{\cond}{{\mbox{\sc CondBm}}}
\newcommand{\bex}{{\mbox{\sc BoxExit}}}
\newcommand{\diff}{{\mbox{\sc DiffExit}}}
\newcommand{\tcu}{{\mbox{\sc CurrentTime}}}
\newcommand{\bdiff}{{\mbox{\sc BanditDiffExit}}}

\newtheorem{thm}{Theorem}[section]

\newtheorem{prop}[thm]{Proposition}

\newtheorem{rem}[thm]{Remark}
\numberwithin{equation}{section}

\title{Exact simulation of diffusion first exit times: algorithm acceleration.}
\begin{document}
\author{S. Herrmann$^1$ and C. Zucca$^2$\\[5pt]
\small{$^1$Institut de Math{\'e}matiques de Bourgogne (IMB) - UMR 5584, CNRS,}\\
\small{Universit{\'e} de Bourgogne Franche-Comt\'e, F-21000 Dijon, France} \\
\small{Samuel.Herrmann@u-bourgogne.fr}\\[5pt]
\small {$^2$Department of Mathematics 'G. Peano', }\\
\small{University of Torino, Via Carlo Alberto 10,
10123 Turin, Italy,}\\ 
\small{cristina.zucca@unito.it}
}
\maketitle
\begin{abstract}
In order to describe or estimate different quantities related to a specific random variable, it is of prime interest to numerically generate such a variate. In specific situations, the exact generation of random variables might be either momentarily unavailable or too expensive in terms of computation time. It therefore needs to be replaced by an approximation procedure. As was previously the case, the ambitious exact simulation of exit times for diffusion processes was unreachable though it concerns many applications in different fields like mathematical finance, neuroscience or reliability. The usual way to describe exit times was to use discretization schemes, that are of course approximation procedures. Recently, Herrmann and Zucca \cite{Herrmann-Zucca-2} proposed a new algorithm, the so-called GDET-algorithm (General Diffusion Exit Time), which permits to simulate exactly the exit time for one-dimensional diffusions. The only drawback of exact simulation methods using an acceptance-rejection sampling is their time consumption. In this paper the authors highlight an acceleration procedure for the GDET-algorithm based on a multi-armed bandit model. The efficiency of this acceleration is pointed out through numerical examples.
\end{abstract}
\par\medskip

 \par\medskip
%
%
%
%
%
%
%
%
%
%
%
%
\section*{Introduction}
A precise description of the first time a given stochastic process exits from a domain is required in many mathematical applications: it can for instance be related to the evaluation of risk of default in mathematical finance or to the description of spike trains in neuroscience,...
Unfortunately, in the diffusion framework (solutions of stochastic differential equations) a simple and explicit expression of the exit time distribution is not attainable except in a few specific cases. It is therefore challenging to find out how to generate such variates. One way to overcome this issue is to introduce an algorithm based on an approximation procedure. Several studies are for instance based on a discretization scheme for the corresponding stochastic differential equation. Most of them are based on improvements of the classical Euler scheme (see for instance \cite{Broadie-Glasserman-Kou-1997}, \cite{Gobet-Menozzi-10}, \cite{Gobet-2000}) which essentially consists in reducing the error stem from the approximation procedure. Another way to deal with the distribution of first exit times consists in approximating their probability density functions and thus in approximating the solution of an integral equation \cite{Sacerdote-2014}. \par\medskip
Apart from all these approximation procedures, Herrmann and Zucca \cite{Herrmann-Zucca-2} proposed an exact simulation of diffusion exit times based on an acceptance-rejection method. The method is directly linked to the Girsanov transformation, a crucial tool already used for the exact simulation of diffusion paths on a fixed time interval \cite{Beskos-2006,beskos2005exact} or for the simulation of first passage times \cite{Herrmann-Zucca}. It is impossible to reasonably compare the numerical methods listed so far since they are of very different types. On the one hand, approximation methods are fast but induce small errors to be controlled. On the other hand, exact method are rather time-consuming.  \par\medskip

The aim of this paper is to improve and accelerate the algorithm presented in \cite{Herrmann-Zucca-2} which permits to generate numerically the first exit time and exit location of a diffusion process from a given interval $[a,b]$. Let us consider the stochastic process $(X_t,\ t\ge 0)$, solution of the SDE:
\begin{equation}\label{eq:depart}
dX_t=\mu(X_t)dt+\sigma(X_t)dB_t,\quad X_0=x\in(a,b),
\end{equation}
where $(B_t,\ t\ge 0)$ stands for the standard one-dimensional Brownian motion, $\sigma\in\mathcal{C}^3([a; b])$ is a positive function on the whole interval $[a,b]$ and $\mu\in\mathcal{C}^2([a; b])$. In the particular case when the function $\sigma$ is constant
we can relax the hypothesis on $\mu$ and we just take $\mu\in \mathcal{C}^1([a; b])$. We denote by $\tau_{a,b}$ the first time the diffusion exits from the interval $[a,b]$:
\begin{equation}\label{eq:def:tau}
\tau_{a,b}(X):=\inf\{t> 0: \ X_t\notin [a,b]  \}.
\end{equation}
Let $T>0$. We call $\bex(x,[a,b],T)$ the efficient algorithm which permits to simulate exactly the random vector $(\tau_{a,b}(X)\wedge T, X_{\tau_{a,b}(X)\wedge T})$, that is the first time the path of the diffusion process $(X_t)_{t\ge 0}$ exits from the time-space rectangle $[0,T]\times [a,b]$ and its associated location. A simple and unified version of this algorithm is presented in Section \ref{sec:Boxexit}, Figure \ref{fig:flowchart} (it corresponds to the algorithms DET and $\kappa$-DET introduced in \cite{Herrmann-Zucca-2}).\medskip

Of course \bex\ is only a basic component for the exit problem from the interval $[a,b]$: the authors suggested in \cite{Herrmann-Zucca-2} to use the Markov property of the time-homogeneous diffusion \eqref{eq:depart} in order to simulate $\tau_{a,b}(X)$. More precisely, the iteration procedure is initialized by $Z_0=x$, the starting position of the diffusion. Then the sequence defined by 
\[
(\mathcal{T}_{n+1}, Z_{n+1})\leftarrow\bex(Z_{n},[a,b],T)
\]
and stopped as soon as $Z_n$ reaches either the value $a$ or $b$ permits to generate the couple $(\tau_{a,b}(X), X_{\tau_{a,b}(X)})$. The efficiency (time consumption) is just related to the unique parameter $T$ since the size of the time-space rectangle associated to the basic component is $[0,T]\times[a,b]$.\medskip 

\begin{figure}[ht]
\centerline{\includegraphics[scale=0.75]{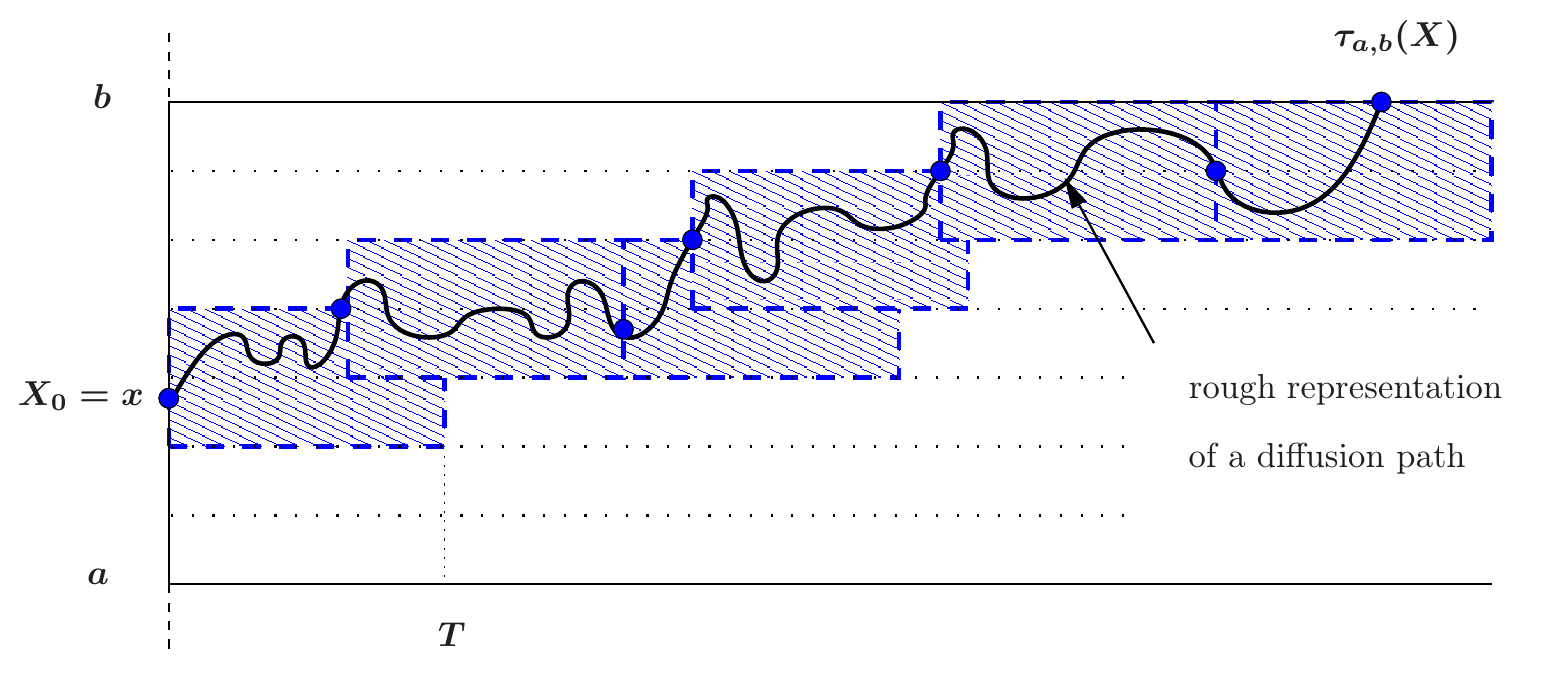}}
\caption{Diffusion path and associated random walk on rectangles}
\label{fig:explanation}
\end{figure}

The main idea of the acceleration procedure is to choose in an optimal way the box size related to the basic components. Instead of fixing the elementary box size equal to $[0,T]\times [a,b]$  ($[a,b]$ being the interval of the initial problem), we propose to cover the interval $]a,b[$ by a fixed number (denoted $N-1$ in Section \ref{sec:randomwalk}) of slices of identical width: 
$]a,b[=\cup_{i=1}^{N-1}I_i$ and to successively use the basic components $\bex(\cdot,I_i,T)$ associated to the family of box sizes $([0,T]\times I_i)_{1\le i\le N-1}$ until the exit of the interval $]a,b[$ occurs. In other words, we introduce a random walk on small rectangles and stop it as soon as it reaches either $a$ or $b$, see Figure \ref{fig:explanation}. At first glance, such a procedure seems to slow down the exact simulation of the exit time since we introduce a new random walk and increase the number of appeals to basic components. But the observation reveals something surprising: for suitable choices of parameters $N$ and $T$, the introduction of the random walk effectively speeds up the algorithm. It is less time-consuming for a diffusion process to exit from boxes of intermediate size compared to boxes of small or large size due to the acceptance-rejection method. This simple argument partly explains the over-performance of the modified algorithm.\medskip
It is therefore challenging to find the optimal parameters $T$ and $N$ in order to obtain the most efficient algorithm. Instead of considering in detail all families of diffusion processes and determining the best choice of parameters on a case-by-case basis, we prefer to propose a randomized algorithmic approach. We find a reasonable value of $T$ and choose $N$ with a multi-armed bandit method ($\epsilon$-greedy algorithm). Such general method can be applied to any diffusion process.\medskip

The material is organized as follows: in Section \ref{sec:Boxexit}, we emphasize a unified and simple version for the exact simulation of exit times, denoted by $\bex$. Section \ref{sec:randomwalk} concerns the introduction of the random walk on small rectangles of area $2T\times [a,b]/N$. A multi-armed bandit method is introduced in Section \ref{sec:multiarmed} for the optimal choice of the parameter $N$. Finally, in the last section we illustrate the efficiency of this new algorithm considering classical diffusion processes like the Ornstein-Uhlenbeck process or the Cox-Ingersoll-Ross model.  
%
%
%
%
%
%
%
%
%
%
%
%
\section{Exit problem from a rectangle}
\label{sec:Boxexit}
Let us first recall the algorithm introduced in \cite{Herrmann-Zucca-2} (see Theorem 4.3) which permits to exactly simulate the exit time from the rectangle $[0,T]\times[l,u]$ for the diffusion path $(X_t,\,t\ge 0)$. The algorithm essentially needs two basic elements:
\begin{enumerate}
\item the exact simulation of the exit time and location $(\mathcal{T},B_\mathcal{T}^x)$ from the interval $[l,u]$ for the Brownian motion starting in $x$. The generation of such a random vector is available (see Section 3 in \cite{Herrmann-Zucca-2}) and will be denoted by $\bre(x,[l,u])$  in the sequel.
\item the generation of the Brownian position $B_t^x$ given $\mathcal{T}>t$ which is denoted $\cond(x,[l,u],t)$ (see  Section 2 in \cite{Herrmann-Zucca-2}).
\end{enumerate}
Both elements allow the construction of a general algorithm for the simulation of exit times. Before introducing the general procedure, we shall focus our attention onto a particular diffusion process which corresponds to the unique solution of a stochastic differential equation with unit diffusion coefficient:
\begin{equation}\label{eq:sigma1}
dX_t=\mu_0(X_t)dt+dB_t,\quad X_0=x\in(a,b).
\end{equation}
Here the drift term is assumed to satisfy $\mu_0\in\mathcal{C}^2([a; b])$.
We define particular functions associated to equation \eqref{eq:sigma1} as:
\begin{align*}
\beta(x)&:=\exp\int_0^x \mu_0(y)\,dy\quad\mbox{and}\quad \gamma(x):=\frac{\mu_0^2(x)+\mu_0'(x)}{2}.
\end{align*}
These functions play an important role in the simulation and do not depend on the considered interval $[l,u]$. Let us now complete these functions with different parameters depending on the interval $[l,u]$: 
\begin{align*}
\beta^+:=\sup_{x\in[l,u]}\beta(x),\quad\gamma^-:=\inf_{x\in[l,u]}\gamma(x)\wedge 0,\quad \gamma^+:=\sup_{x\in[l,u]}\gamma(x), \quad \gamma^0:=\gamma^+-\gamma^-.
\end{align*}
 
A unified statement of the exact simulation algorithms presented in \cite{Herrmann-Zucca-2} is defined as follows: 
\begin{prop} The couple $(\tau_{l,u}(X)\wedge T, X_{\tau_{l,u}(X)\wedge T})$ which corresponds to the exit problem of the diffusion path \eqref{eq:sigma1} from the rectangle $[0,T]\times[l,u]$, has the same distribution than the outcome $(\mathcal{T},Z)$ of the algorithm $\bex(x,[l,u],T)$ for any $T>0$ (see the flowchart in Figure \ref{fig:flowchart}).
\end{prop}
It is worth noting that the random variables generated in the algorithm $\bex$ (i.e. $E$, $U$, $V$, $W$) are independent (In Figure \ref{fig:flowchart}, $U$, $V$ and $W$ are represented by the same character $U_\bullet$ which corresponds to independent uniformly distributed variates). 
\begin{figure}
\centering
\centerline{\tikz \node [scale=0.6, inner sep=0] {\begin{tikzpicture}[node distance=1.7cm, every node/.style={fill=white, font=\sffamily}, align=center]
  \node (A1) [activityStarts] {initialization\\ $\quad\mathcal{T}=0$, $Z=x$, $K=T\quad$};
  \node (A2) [process, left of=A1,yshift=-2cm] {generate\\ $(S,Y)= \bre(Z,[l,u])$\\[4pt]
  \ generate $E\sim\mathcal{E}(\gamma^0)$};
  \node (A3) [startstop, below of=A2, yshift=-0.2cm] {test $S=\min(K,E,S)$};
  \node (A4) [startstop, below of=A3,yshift=-0.7cm] {test $\beta^+ U_\bullet\le \beta(Y)$ and\\[4pt] $\log(U_\bullet)\le \gamma^-\,(K-S)$};
  \node (A5) [activityRuns, below of=A4, yshift=-0.7cm] {$\mathcal{T}\leftarrow \mathcal{T}+S$\\ $Z=Y$};
  \node (A6) [activityStarts, right of=A5, xshift=8.5cm] {\ \ outcome: $\mathcal{T}$ and $Z$ };
  \node (B1) [startstop, right of=A2, xshift=3.5cm] {test $K=\min(K,E,S)$};
  \node (B2) [process, below of=B1, yshift=-0.5cm] {generate\\ $Y_c=\cond(Z,[l,u],K)$};
  \node (B3) [startstop, below of=B2] {test $\beta^+ U_\bullet\le \beta(Y_c)$};
  \node (B4) [activityRuns, below of=B3,yshift=-0.5cm] {$\mathcal{T}\leftarrow \mathcal{T}+K$\\ $Z=Y_c$};
  \node (C1) [process, right of=B3, xshift=4.5cm] {generate\\ $Y_c=\cond(Z,[l,u],E)$};
  \node (C2) [startstop, above of=C1] {test $\gamma^0 U_\bullet> \gamma (Y_c)-\gamma_-$};
  \node (C3) [activityRuns, above of=C2,yshift=0.2cm] {$\mathcal{T}\leftarrow \mathcal{T}+E$\\ $Z=Y_c$,\ $K\leftarrow K-E$};
  \draw[->] (A1) -- (A2);
  \draw[->] (B1) --++ (3,0) --++ (0,-3.5) -- node [xshift=-0.3cm,yshift=2.1cm]{No}(C1);
  \draw[->] (A2) -- (A3);
  \draw[->] (A3) -- node[text width=0.5cm]{Yes} (A4);
  \draw[->] (A3)  --++ (2,0.5) -- node [xshift=-0.2cm,yshift=-0.1cm]{No} (B1);
  \draw[->] (A4) -- node[text width=0.5cm]{Yes} (A5);
  \draw[->] (A5) -- ++(2,-0.2)-- ++(5.5,0) -- (A6);
  \draw[->] (B1) -- node[text width=0.5cm]{Yes} (B2);
  \draw[->] (B2) -- (B3);
  \draw[->] (B3) -- node[text width=0.5cm]{Yes} (B4);
  \draw[<-] (A1) -- ++(-4.7,0) -- ++(0,-6.3) -- node{No}(A4);    
  \draw[->] (B4) -- (A6);
  \draw[->] (C1) -- (C2);
  \draw[->] (C2)  --++ (2.5,0) --++ (0,4.1)-- node [xshift=-0.2cm,yshift=-0cm]{No} (A1); 
  \draw[->] (C2) -- node[text width=0.5cm]{Yes} (C3);
  \draw[->] (C3)  --++ (0,1.2)  --++ (-9.5,0) -- (A2);
  \draw[->] (B3)  --++ (2,-1.2)  --++ (7.5,0) --++ (0,7.5) --++ (-9.6,0) -- node[xshift=5.2cm,yshift=-7.5cm]{No} (A1);                          
  \end{tikzpicture}};}
\caption{Flowchart of the algorithm  $\bex(x,[l,u],T)$} \label{fig:flowchart}
\end{figure}

{\small
\begin{oframed}
\begin{algorithm}[H]
\SetAlgoLined
 \KwData{ $x$ (starting position), $T$, $l$ and $u$ (box size),  $\gamma(\cdot)$ and $\beta(\cdot)$ (input functions).}
 \KwResult{the random time $\mathcal{T}$ and the random location $Z$.}
\vspace*{0.2cm} 
Initialization: $K=T$, $Z=x$, $\mathcal{T}=0$, ${\rm test}=0$\;
Computation of $\gamma^-$, $\gamma^0$, $\beta^+$ depending on the interval $[l,u]$\;
\vspace*{0.2cm} 
 \While{${\rm test}=0$}{
 generate $E\sim\mathcal{E}(\gamma^0)$ and $U\sim V\sim W\sim \mathcal{U}([0,1])$\;
 generate $(S,Y)= \bre(Z,[l,u])$\;
  \uIf{$S=\min(K,E,S)$}{
  \eIf{$\beta^+ U\le \beta(Y)$ {\rm and} $\log(W)\le \gamma^-\,(K-S)$}{set ${\rm test}=1$, $Z\leftarrow Y$ and $\mathcal{T}\leftarrow \mathcal{T}+S$\;}{go to \emph{initialization}\; }
   }
   \uElseIf{$K=\min(K,E,S)$}{generate
\(
Y_c=\cond(Z,[l,u],K)
\)\;
\eIf{$\beta^+ U\le \beta(Y_c)$}{set ${\rm test}=1$, $Z\leftarrow Y$ and $\mathcal{T}\leftarrow \mathcal{T}+K$\;}{go to \emph{initialization}\;}}
   \Else{generate
\(
Y_c=\cond(Z,[l,u],E)
\)\;
\eIf{$\gamma^0 V> \gamma (Y_c)-\gamma_-$}{$Z\leftarrow Y_c$, $\mathcal{T}\leftarrow \mathcal{T}+E$ and $ K\leftarrow K-E$\;}{go to \emph{initialization}\;} 
 }}
 \caption{$\bex(x,[l,u],T)$}
\end{algorithm}
\end{oframed}}
\begin{rem}
\label{rem:Tinfinite} Under the assumption $\gamma^-=0$ that is $\inf_{x\in[l,u]}\gamma(x)\ge 0$, it is allowed to choose $T=\infty$ in the algorithm $\bex(x,[l,u],T)$. It sould be noted that $\bex$ with $T<\infty$ corresponds to the so-called \emph{$\kappa$-DET} algorithm in \cite{Herrmann-Zucca-2} whereas $\bex$ with $T=\infty$ corresponds to the \emph{DET} algorithm. Here we decided to unify the presentation for pedagogical reasons.
\end{rem}
\begin{rem} 
\label{rem:lamperti}
The Lamperti transform permits to generalize the study to equations with non-unitary diffusion coefficients as \eqref{eq:depart}. 
We simply present this well-known transformation. Let $(X_t)_{t\ge 0}$ be the unique solution to the SDE \eqref{eq:depart} and let us introduce 
\begin{equation}
\label{eq:defofS}
\mathcal{S}(x)=\int_0^x\frac{du}{\sigma(u)},\quad \forall x\in\mathbb{R},
\end{equation}
then It\^o's lemma implies that $\widehat{X}_t:=\mathcal{S}(X_t)$ satisfies \eqref{eq:sigma1} with initial condition  $\widehat{X}_0=\mathcal{S}(X_0)$ and drift term
\[
\mu_0(x):=\frac{\mu(\mathcal{S}^{-1}(x))}{\sigma(\mathcal{S}^{-1}(x))}-\frac{1}{2}\sigma'(\mathcal{S}^{-1}(x)),\quad x\in\mathbb{R}.
\]
The procedure to simulate the exit time and location of a diffusion path $(X_t)_{t\ge 0}$ defined by \eqref{eq:depart} from the rectangle $[0,T]\times[l,u]$ is therefore the following:
\begin{enumerate}
\item Simulate $(\mathcal{T},Z)$ the exit time and location of the diffusion $(\widehat{X}_t)_{t\ge 0}$ using the algorithm $\bex(\mathcal{S}(x),[\mathcal{S}(l),\mathcal{S}(u)],T)$
\item Compute $\mathcal{S}^{-1}(Z)$. Then $(\mathcal{T}, \mathcal{S}^{-1}(Z))$ corresponds to the exit time and location of the diffusion $(X_t)_{t\ge 0}$ from the interval $[l,u]$.
\end{enumerate}
\end{rem}
%
%
%
%
%
%
%
%
%
%
%
%
\section{A random walk on rectangles} 
\label{sec:randomwalk}
Using the exit problem of rectangles as the basic component, we can build a general algorithm that enables us to simulate exactly the exit time of the diffusion process \eqref{eq:depart} from the interval $[a,b]$. Applying the Lamperti transformation already described in Remark \ref{rem:lamperti}, there is a one-to-one correspondence between the process $(X_t)$ solution of \eqref{eq:depart} starting in $X_0=x$ and $(\widehat{X}_t)$ the solution of \eqref{eq:sigma1} starting in $\widehat{X}_0=\mathcal{S}(x)=\hat{x}$ where $\mathcal{S}$ is defined by \eqref{eq:defofS}. Moreover, the interval $[a,b]$ is transformed into $[\hat{a},\hat{b}]=[\mathcal{S}(a), \mathcal{S}(b)]$.\medskip

Let us now describe how to deal with the exit problem for $(\widehat{X}_t)$ associated to the interval $[\hat{a},\hat{b}]$. Let us first fix a parameter $T>0$ and a number $N\ge 2$ (we shall comment on these choices later on). These parameters define the size of the typical boxes used in the algorithm illustrated by Figure \ref{fig:explanation}: rectangles of area $2(\hat{b}-\hat{a})T/N$. The main idea is quite simple: the interval $[\hat{a},\hat{b}]$ is split into $N$ intervals of identical length $\delta$, associated to the following space grid: $a_0=\hat{a}$ and $a_{j+1}=a_j+\delta$ for $0\le j\le N-1$. Here $\delta=(\hat{b}-\hat{a})/N$. We define the index function: 
\begin{equation}
\label{eq:index}
\imath(x)=j \quad \mbox{if}\quad  (x-\hat{a})\in\left[\frac{\delta}{2}+(j-1)\delta,\frac{\delta}{2}+j\,\delta \right[,
\end{equation}
otherwise either $\imath(x)=1$ for $x\le \hat{a}+\delta/2$ or $\imath(x)=N-1$ for $x\ge \hat{b}-\delta/2$.\\[2pt]
\centerline{\begin{tikzpicture}[y=1cm, x=1cm, thick, font=\footnotesize]    
\usetikzlibrary{arrows,decorations.pathreplacing}
\tikzset{
   brace_top/.style={
     decoration={brace},
     decorate
   },
   brace_bottom/.style={
     decoration={brace, mirror},
     decorate
   }
}
\draw[line width=1.2pt, >=latex'](0,0) -- coordinate (x axis) (8,0) node[right] {}; 
\draw (0,0.1) -- (0,-0.1) node[below,pos=1.5] {$\hat{a}$};
\draw (1,0.1) -- (1,-0.1) node[below] {$\hat{a}+\delta$};
\draw (2,0.1) -- (2,-0.1) node[below] {$\hat{a}+2\delta$};
\draw (3,0.1) -- (3,-0.1) node[below] {$\hat{a}+3\delta$};
\draw (4,0.1) -- (4,-0.1) node[below,pos=1.5] {$\ldots$};
\draw (6,0.1) -- (6,-0.1) node[below,pos=1.5] {$\ldots$};
\draw (7,0.1) -- (7,-0.1) node[below] {$\hat{b}-\delta$};
\draw (8,0.1) -- (8,-0.1) node[below] {$\hat{b}$};
\node (start_week) at (0,0.1) {};
\node (end_week) at (1.5,0.1) {};
\draw [brace_top] (start_week.north) -- node [above, pos=0.5] {$\imath=1$} (end_week.north);
\node (start_week1) at (1.5,0.1) {};
\node (end_week1) at (2.5,0.1) {};
\draw [brace_top] (start_week1.north) -- node [above, pos=0.5] {$\imath=2$} (end_week1.north);
\node (start_week2) at (2.5,0.1) {};
\node (end_week2) at (3.5,0.1) {};
\draw [brace_top] (start_week2.north) -- node [above, pos=0.5] {$\imath=3$} (end_week2.north);
\node (start_week3) at (5.5,0.1) {};
\node (end_week3) at (6.5,0.1) {};
\draw [brace_top] (start_week3.north) -- node [above, pos=0.5] {$\imath=N-2\ \ $} (end_week3.north);
\node (start_week4) at (6.5,0.1) {};
\node (end_week4) at (8,0.1) {};
\draw [brace_top] (start_week4.north) -- node [above, pos=0.5] {$\ \ \ \ \imath=N-1$} (end_week4.north);
\end{tikzpicture}}
Each index value $\imath\in \{1,2,\dots,N-1\}$ is associate to an interval of length $2\delta$: 
\begin{equation}\label{eq:interval}
I_{\imath}=]\hat{a}+(\imath-1)\delta,\hat{a}+(\imath+1)\delta[.
\end{equation} 
We notice that the family of intervals $(I_\imath)_{1\le \imath\le N-1}$ is a covering of the initial interval $]\hat{a},\hat{b}[$. Moreover, for any $x\in]\hat{a},\hat{b}[$, $x\in I_{\imath(x)}$. \medskip

A random walk  corresponding to a skeleton of the diffusion path can be thus constructed (see Figure \ref{fig:explanation}): $(T_0,Y_0)=(0,\hat{x})$ is the starting time and position of the diffusion process $(\widehat{X}_t)$, solution of \eqref{eq:sigma1}. The random sequence $(T_{n+1},Y_{n+1})$ is defined recursively as follows: $T_{n+1}-T_n$ stands for the exit time of the diffusion starting in $Y_{n}$ from the rectangle $[0,T]\times I_{\imath(Y_n)}$ and $Y_{n+1}$ corresponds to the associated exit location. Let us  define \[\mathcal{N}:=\inf\{n\ge 0:\ Y_n\notin ]\hat{a},\hat{b}[\}\] then the combination of the Markov property and the Lamperti transform implies the following statement.
\begin{prop}\label{prop2}
The diffusion exit time and location $(\tau_{a,b}(X),X_{\tau_{a,b}(X)})$ has the same distribution as the stopped random walk $(T_{\mathcal{N}},\mathcal{S}^{-1}(Y_{\mathcal{N}}))$ and consequently the same distribution as $(\mathcal{T},Z)$ the outcome of the algorithm $\diff$. 
\end{prop} The algorithm $\diff$ induced by this random walk is the following.
\begin{framed}
\begin{algorithm}[H]
\SetAlgoLined
 \KwData{ $x$ (starting position of the diffusion), $T$, $N$ (box size), $\gamma(\cdot)$ and $\beta(\cdot)$ (input functions), $\mathcal{S}(\cdot)$ (Lamperti transform).}
 \KwResult{the random time $\mathcal{T}$ and the random location $Z$.}
\vspace*{0.2cm} 
initialization: $\mathcal{T}=0$, $Z=\mathcal{S}(x)$, $\hat{a}=\mathcal{S}(a)$, $\hat{b}=\mathcal{S}(b)$\;
 \While{$Z\in]\hat{a},\hat{b}[$}{$(S,Z)\leftarrow\bex(Z,I_{\imath(Z)},T)$\;
 $\mathcal{T}\leftarrow\mathcal{T}+S$\;
 }
 $Z\leftarrow \mathcal{S}^{-1}(Z)$\;
 \caption{Diffusion Exit Problem $\diff(T,N)$}
 \label{algo:2}
\end{algorithm}
\end{framed} 
Of course, the efficiency of this exact simulation algorithm heavily depends on the parameters $T$ and $N$ which characterize the size of the typical boxes. If the box is large, then the algorithm $\bex$ becomes time consuming since it is based on a rejection sampling. On the contrary, small boxes imply that the random walk on rectangles requires a lot of iterations in order to hit the boundaries of the interval $[a,b]$. There is therefore an intermediate box size which permits to observe simulations that take a reasonable computation time.  

\begin{figure}[h]
\centerline{\includegraphics[width=6cm]{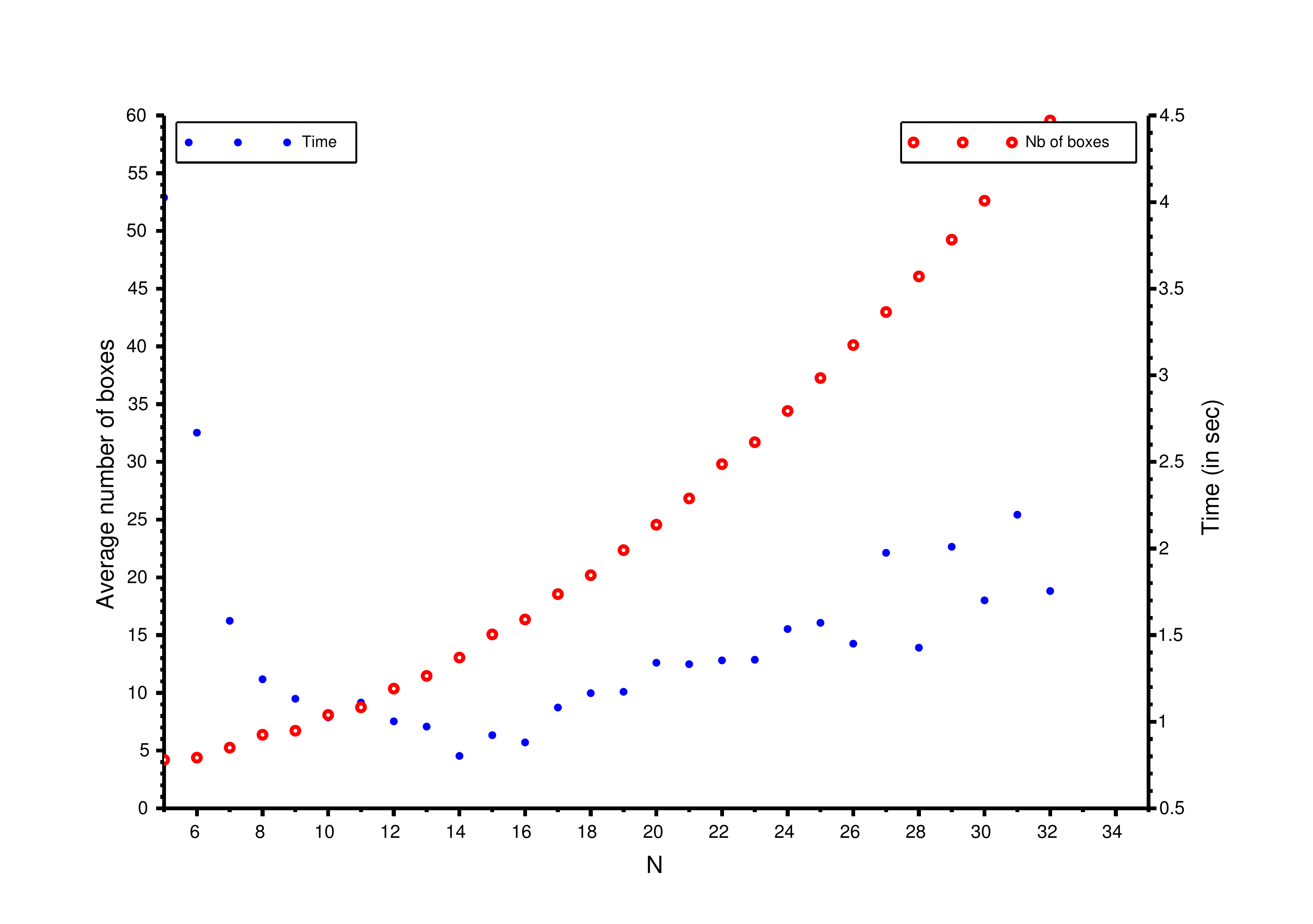}}
\caption{\small Average number of boxes used in the exit algorithm and total computation time (for the simulation of the whole sample) versus the box size parameter $N$ for the diffusion process of Example 1.  Exit problem from the interval $[a,b]=[0,7]$ with the starting position $x=3$. Each value is obtained with a sample of size $10\,000$ and $T=1$.}
\label{Fig1}
\end{figure}

In order to illustrate this feature, let us introduce two particular examples:\\ 
{\bf Example 1:} the diffusion process with unitary diffusion coefficient and with the following drift term: $\mu_0(x)=2+\sin(x)$. We consider the exit problem from the interval $[a,b]=[0,7]$, the diffusion starting in $x=3$. Figure \ref{Fig1} represents on one hand the average number of boxes needed in order to observe the exit depending on the box size (we let $N$ vary). On the other hand we also point out the computation time (in sec) needed to generate a sample of 10\,000 
diffusion exit times.\\

\begin{figure}[h]
\centerline{\includegraphics[width=6cm]{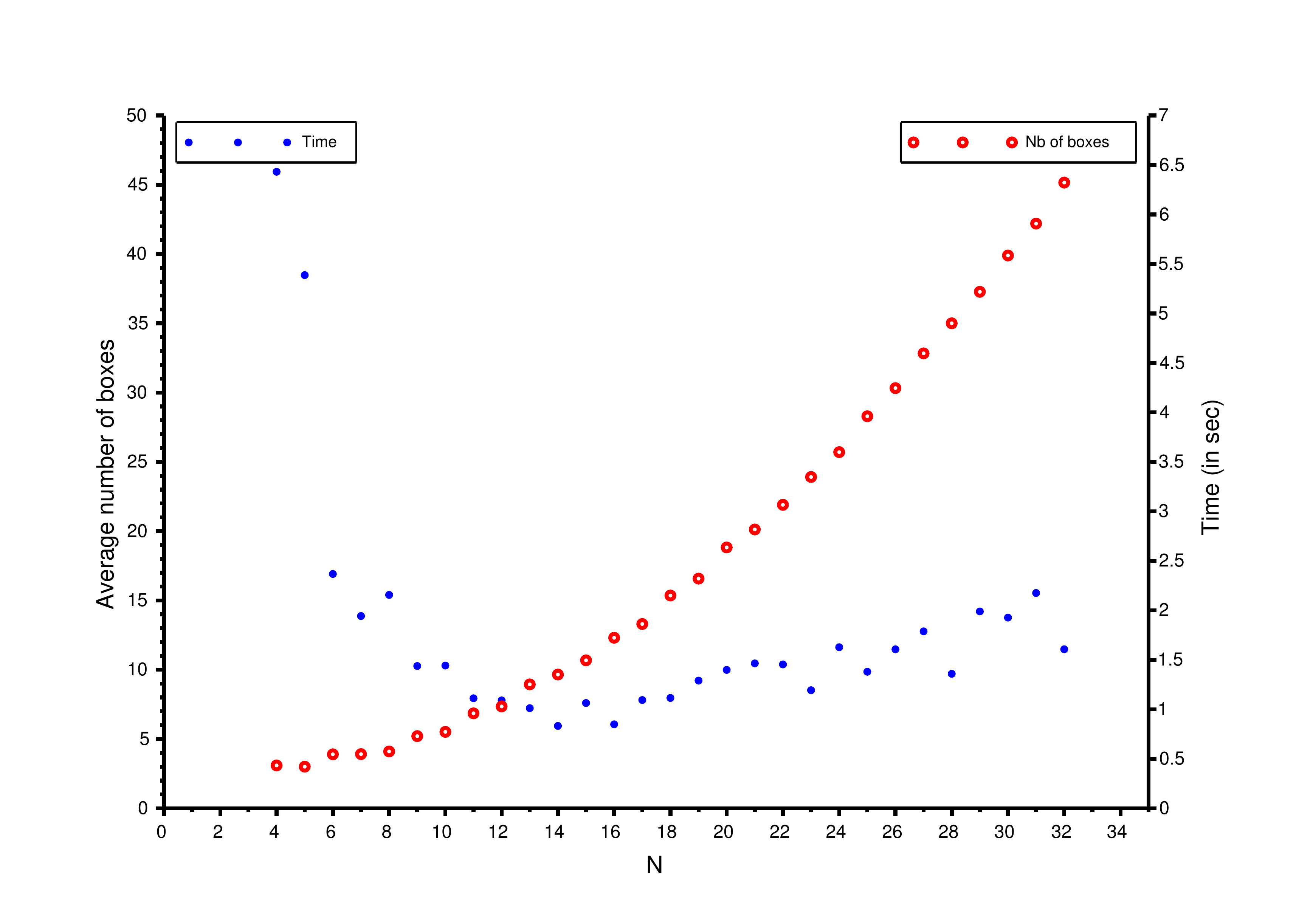}\hspace*{-0.5cm} \includegraphics[width=6cm]{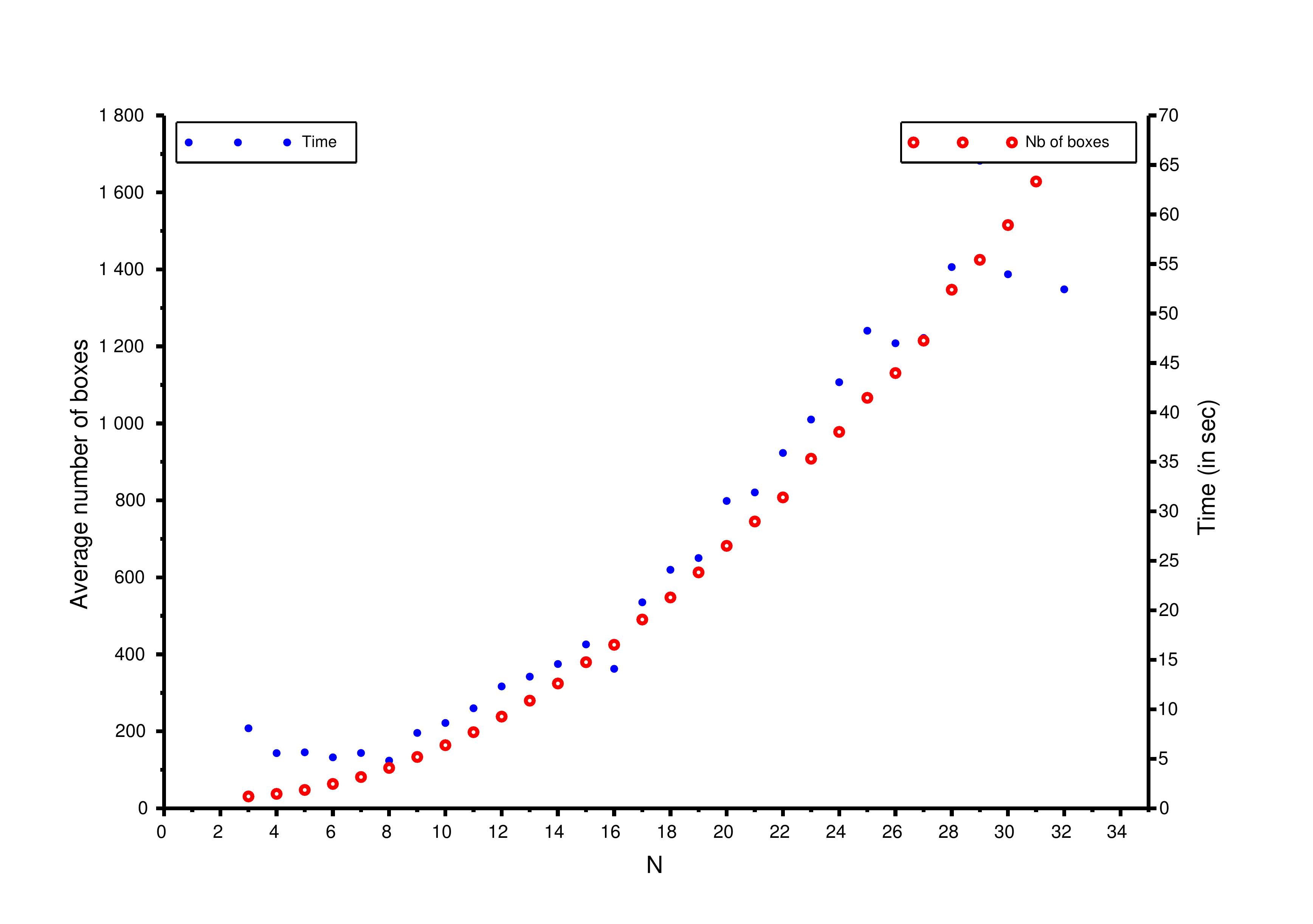} }
\caption{\small Average number of boxes and total computation time versus the box size parameter $N$ for the Ornstein-Uhlenbeck process with parameter $\lambda=1$. Exit problem from the interval $[a,b]=[0,7]$ with the starting position $x=3$ (left) and $[a,b]=[-2,2]$ and $x=0.5$ (right). Each value is obtained with a sample of size $10\,000$ and $T=1$.}
\label{Fig2}
\end{figure}

\noindent {\bf Example 2:} the Ornstein-Uhlenbeck process with unitary diffusion coefficient and drift term : $\mu_0(x)=-\lambda x$ with $\lambda>0$. First we focus our attention to the exit problem from the interval $[a,b]=[0,7]$ with the initial condition $x=3$ and the parameter $\lambda=1$, see Figure \ref{Fig2} (left). 

We notice that the optimal box size corresponds to $N=14$ when $T=1$ is fixed. Such an optimal choice strongly depends on the interval $[a,b]$.  Since the diffusion is mean-reverting, let us observe what happens when the interval $[a,b]$ contains $0$. Figure \ref{Fig2} (right) illustrates that $N=5$ is optimal for $[a,b]=[-2,2]$ and $x=0.5$. We also notice that the number of boxes used in such a particular situation is much larger than in the previous situation. It is therefore difficult to obtain a theoretical optimal value for the parameter $N$. That is why we aim to find an acceleration method for the simulation of exit times (Algorithm \ref{algo:2}: $\diff$) using an algorithmic approach based on a multi-armed bandit.
%
%
%
%
%
%
%
%
%
%
%
%
\section{Algorithm acceleration: a multi-armed bandit approach}
\label{sec:multiarmed}
Let us now  suggest an acceleration method for the algorithm $\diff$ presented in the previous section and depending on both parameters $N$ and $T$ (size of the typical boxes). The procedure is the following: we first fix $T>0$ and $N_0\ge 2$. Then we introduce an algorithm used for the multi-armed bandit problem in order to choose an interesting value of $N$ satisfying $N\le N_0$ and reducing the time consumption of the algorithm $\diff$. \medskip

The multi-armed bandit corresponds to a famous problem where reinforcement learning plays a crucial role, theoretical and practical studies aim to find trade-offs between exploration and exploitation. The historical problem is quite simple and related to a slot machine with a finite number of levers. One is faced repeatedly to a choice between these actions and after each choice one receive a random numerical reward depending on the selected lever. The objective is to maximize the average cumulative reward of a series of actions (for instance, $10\,000$ successive selections) using a strategy based on an exploration-exploitation algorithm. The exploration consists in selecting several times any arm of the bandit in order to estimate the different mean rewards while the exploitation focuses on the choice of the arm whose estimated reward is maximal. We refer to the interesting textbooks \cite{Slivkins} and \cite{Sutton-Barto} for practical and theoretical results associated with this reinforcement learning framework. Several bandit algorithms permit to obtain theoretical bounds of the total expected regret which represents a simple performance measure in such a framework: $\epsilon$-greedy, Boltzmann exploration, UCB (Upper Confidence Bounds), etc. Here we  focus our attention on the \emph{$\epsilon$-greedy algorithm} which is rather intuitive, simple to implement and outperforms theoretically
sound algorithms on most settings \cite{Vermorel-Mohri}.\medskip

In our particular situation, the multi-armed bandit corresponds to the algorithm $\diff(T,N)$: each arm represents a value of $N\in \{2,3,\ldots,N_0\}$ which characterizes the space splitting used in the algorithm.  The reward associated with each arm is the numerical time consumption of each exit time generation. It is of course random since the basic components of the algorithm use rejection sampling. Let us mention that the objective is here opposite: the aim is to minimize the cumulative reward... That means that each use of the algorithm $\diff$ leads to an evaluation of the time spent. We shall therefore use a clock for determining the current time denoted by $\tcu$. \medskip

Let us present the application of $\epsilon$-greedy algorithm in such a context. After the $n$-th use of the algorithm $\diff(T,N)$, the empirical mean of the time consumption is denoted by $\bar{\mu}_n(N)$, for $2\le N\le N_0$, and the number of times we already used the arm $N$ until $n$ is $m_n(N)$. In the $\epsilon$-greedy algorithm, the choice of the parameter $N$ evolves randomly as the number of simulations increases and depends on a fixed parameter $\epsilon$. The probability to choose the arm $N$ for the $n$-th simulation is defined by:
\begin{equation}
\label{eq:bandit}
\pi_{n+1}(N)=(1-\epsilon)\mathds{1}_{ \displaystyle \{N=\underset{2\le j\le N_0}{\operatorname{arg\,min}}\ \bar{\mu}_n(j) \}}+\frac{\epsilon}{N_0-1},
\end{equation}
with the starting values $\pi_1(N)=\epsilon/(N_0-1)$ for all $N\in\{2,3,\ldots,N_0\}$. 
Such strategy for the random choice of the parameters permits to globally reduce the consumption time for a sequential use of the algorithm $\diff$. Of course the parameter $\epsilon$ characterizing the competition between exploration and exploitation has an influence on the acceleration strength and should depend on the sample size. Different studies even suggest to let $\epsilon$ depend on the number of actions $\epsilon:=\epsilon(n)$ of the order $\epsilon(n)= n^{-1/3}((N_0-1)\log(n))^{1/3}$ (see, for instance, Theorem 1.4 in \cite{Slivkins}). Nevertheless experimental results emphasize that making the $\epsilon$ decrease does not significantly improve the performance of the multi-armed bandit strategy \cite{Vermorel-Mohri}. In the following we shall therefore only use $\epsilon$-greedy algorithm with fixed value for $\epsilon$.

The modification of the $\diff$ leads to the following algorithm.
\begin{framed}
\begin{algorithm}[H]
\SetAlgoLined
 \KwData{ $x$ (starting position), $T$, $N_0$, $\gamma(\cdot)$ and $\beta(\cdot)$ (input functions), $M$ (size of the sample: number of simulations).}
 \KwResult{Sample of $M$ simulations for the couple random time $\mathcal{T}$ and random location $Z$.}
initialization:\\ $\pi(N)\leftarrow 1/(N_0-1)$, $\bar{\mu}(N)\leftarrow 0$ and $m(N)=0$ for all $2\le N\le N_0$\;
\For{$j\leftarrow 1$ \KwTo $M$}{
choose randomly $N$ w.r.t. the distribution $\pi(\cdot)$\;
$t\leftarrow \tcu$\;
$(\mathcal{T}_j,Z_j)\leftarrow\diff(T,N)$\;
$t\leftarrow \tcu-t$\;
$\bar{\mu}(N)\leftarrow (m(N)\bar{\mu}(N)+t)/(m(N)+1)$\;
$m(N)\leftarrow m(N)+1$\;
\For{$i\leftarrow 2$ \KwTo $N_0$}{
$\pi(i)\leftarrow \epsilon/(N_0-1)$\;}
$\pi(\arg\min \bar{\mu})\leftarrow \pi(\arg\min \bar{\mu})+(1-\epsilon)$\;
}
 \caption{ $\bdiff(T,N_0)$}
 \label{algo:2}
\end{algorithm}
\end{framed} 
This new algorithm called $\bdiff$ outperforms the exact algorithms introduced for the simulation of diffusion exit times in \cite{Herrmann-Zucca-2} as it appears obvious in the numerical illustrations presented in Section \ref{sec:illu}.
%
%
%
%
%
%
%
%
%
%
%
%
%
%
%
%
%
%
%
%
%
%
%
%
\section{Numerical illustration} 
\label{sec:illu}
\subsection{First example}\label{sec:ex1}  First we consider the exit time from the interval $[a,b]$ for the diffusion:
\begin{equation}\label{eq:ex1}
dX_t=(2+\sin(X_t))\,dt+dB_t,\quad t\ge 0, \quad X_0=x.
\end{equation}
In \cite{Herrmann-Zucca-2}, the DET-algorithm permits to generate the exit time due to an acceptance rejection procedure (this algorithm corresponds to the already presented $\bex(x,[a,b],T)$ for the particular value $T=\infty$, we can observe that the condition described in Remark \ref{rem:Tinfinite} is satisfied). Using a sample of exit time generations we can estimate the average computation time. Here the data correspond to the exit time from the interval $[0,7]$ when starting in $x=3$.\\[5pt]
\renewcommand{\arraystretch}{1.2}
\centerline{\begin{tabular}{|c||p{3cm}||p{2cm}|p{2cm}|}
\hline
sample size & average time (ms) &  \multicolumn{2}{|c|}{confidence interval ($95\%$)}\\
\hline
$10\,000$ & $7.832$ & $7.676$ & $7.989$ \\
\hline
\end{tabular}}
\vspace*{0.2cm}

\begin{figure}[h]
\centerline{\includegraphics[width=8cm]{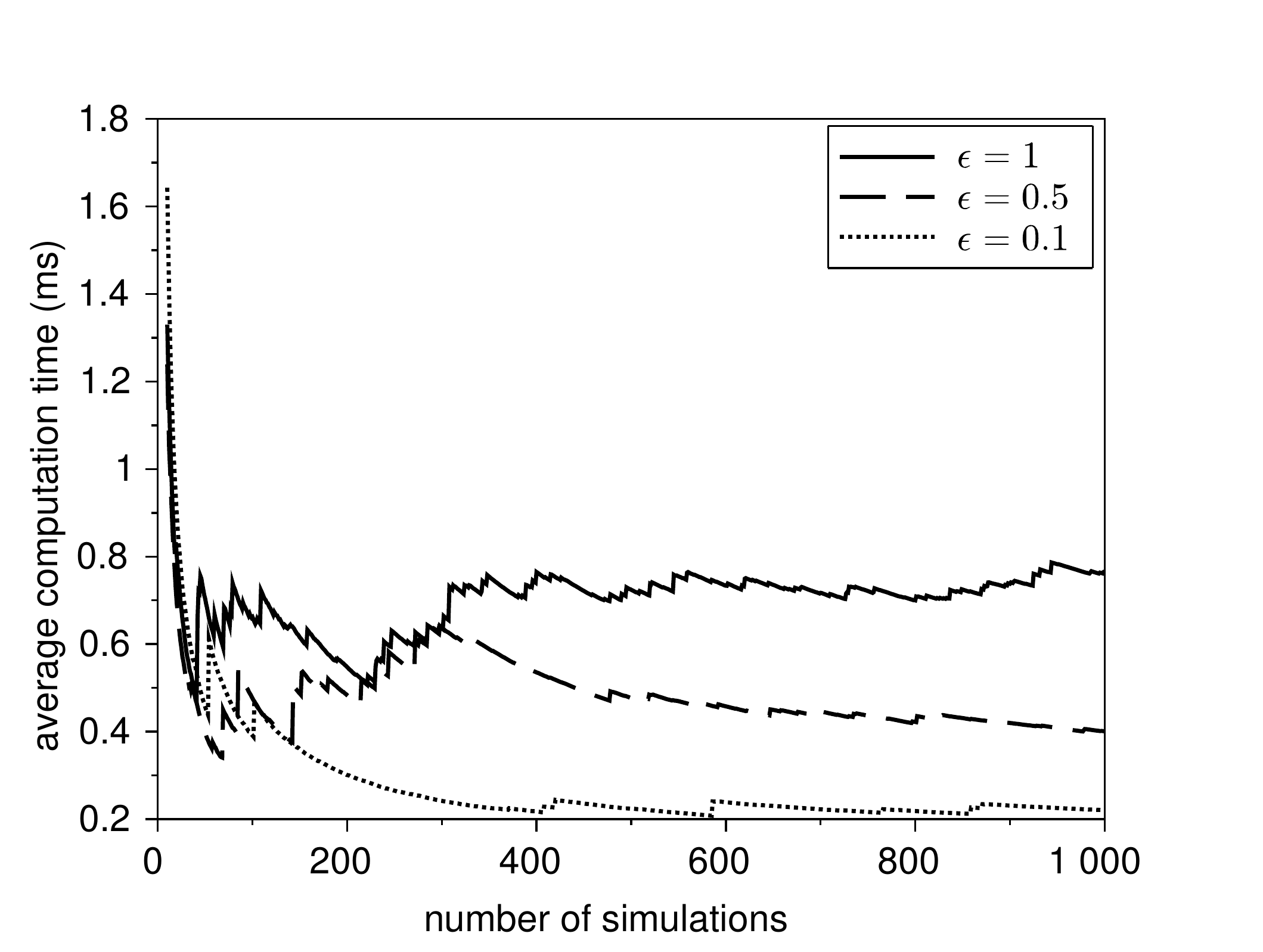}}
\caption{\small Empirical mean of the computation times (in \emph{ms}) versus the number of simulations ($10$ to $1\,000$). Each computation time concerns the simulation of an exit from the interval $[a,b]=[0,7]$ with starting value $x=3$. We use the $\epsilon$-greedy bandit algorithm with different values $\epsilon$ ($\epsilon=1$  corresponds to a uniform choice of the parameter $N$ in $\{2,\ldots, 21\}$). The elementary box size is $2(b-a)T/N$ with $T=1$.  }
\label{Fig1n}
\end{figure}

It is of prime interest to compare the computation time using $\bex$-algorithm with the computation time using the bandit algorithm presented in Section \ref{sec:multiarmed}. Here we deal with a sample of $1\,000$ actions in the bandit algorithm, each run corresponds to the simulation of an exit time from the interval $[a,b]=[0,7]$. Let us note that inbetween two consecutive runs, the bandit algorithm proceed to an optimisation computation corresponding to the choice of the box size. Therefore the sequence of the consumption times $\tau^{(1)},\ldots,\tau^{(n)}$ do not represent i.i.d random variables (the confidence interval is therefore not available). We point out the performance of such an algorithm in Figure \ref{Fig1n}: the averaged computation time is strongly reduced. The figure represents the sequence : $(\frac{1}{n}\sum_{i=1}^n\tau^{(i)})_{n}$ for $10\le n\le 1\,000$. 

\begin{figure}[h]
\centerline{\includegraphics[width=6cm]{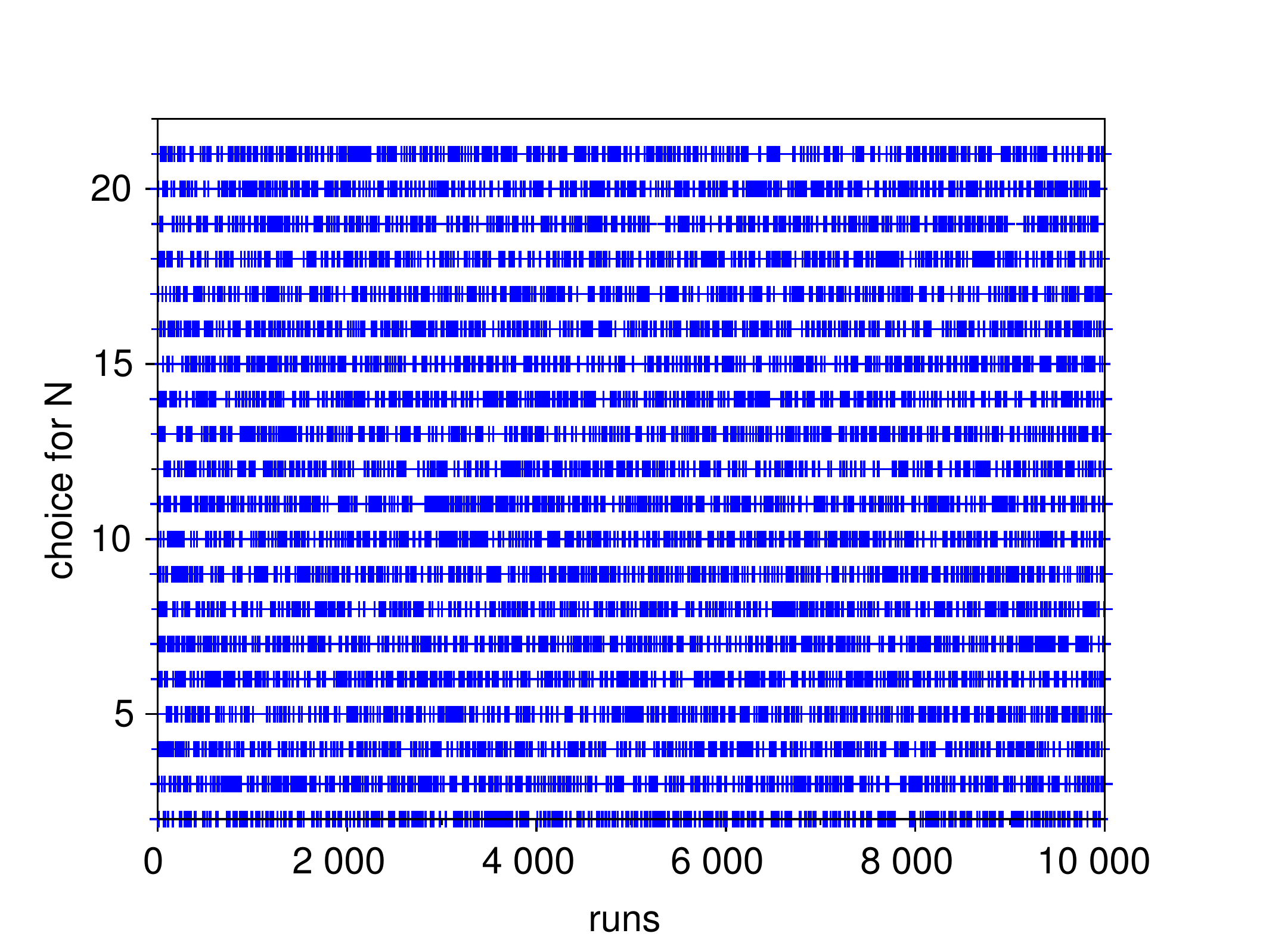}\includegraphics[width=7cm]{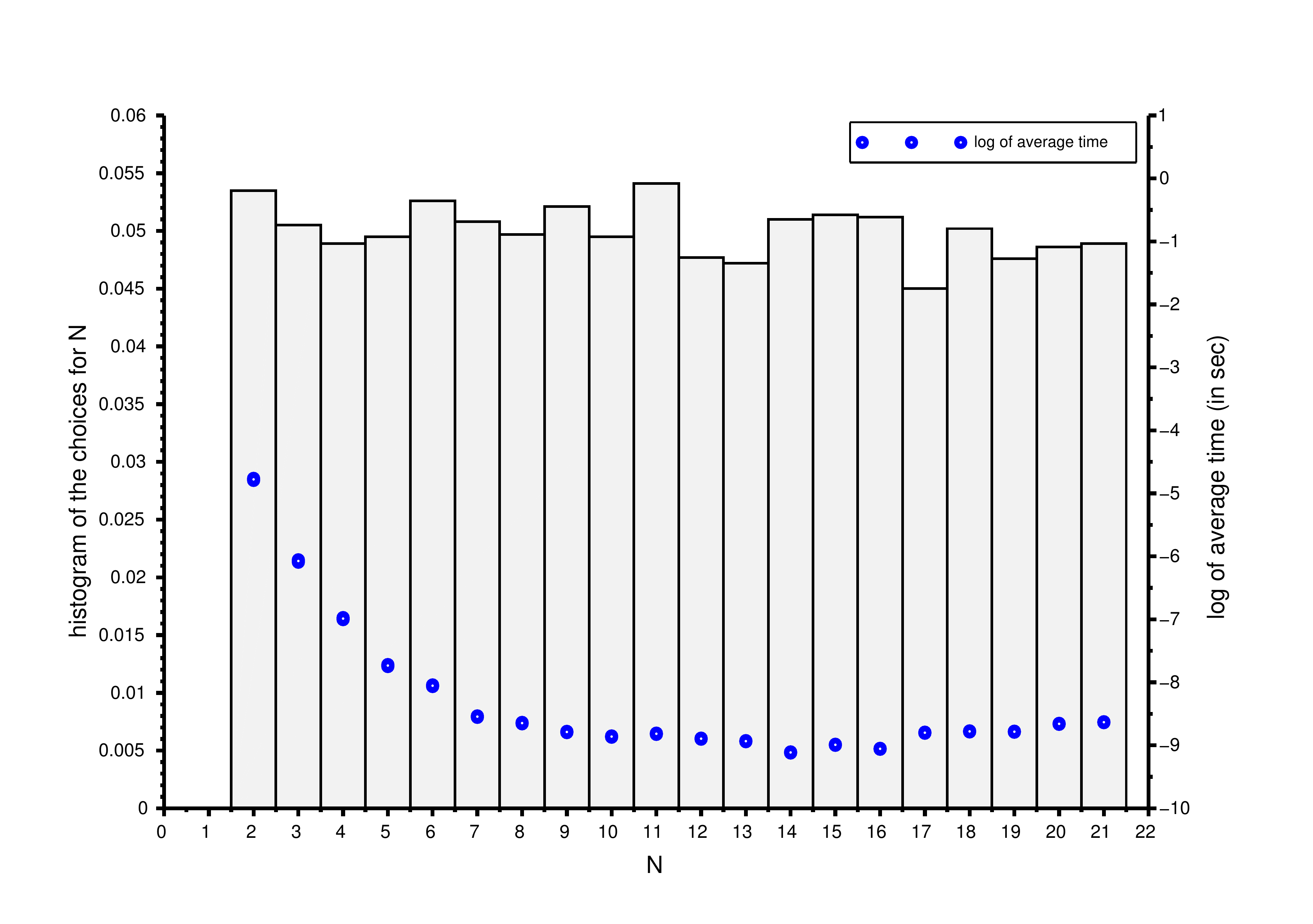}}
\centerline{\includegraphics[width=6cm]{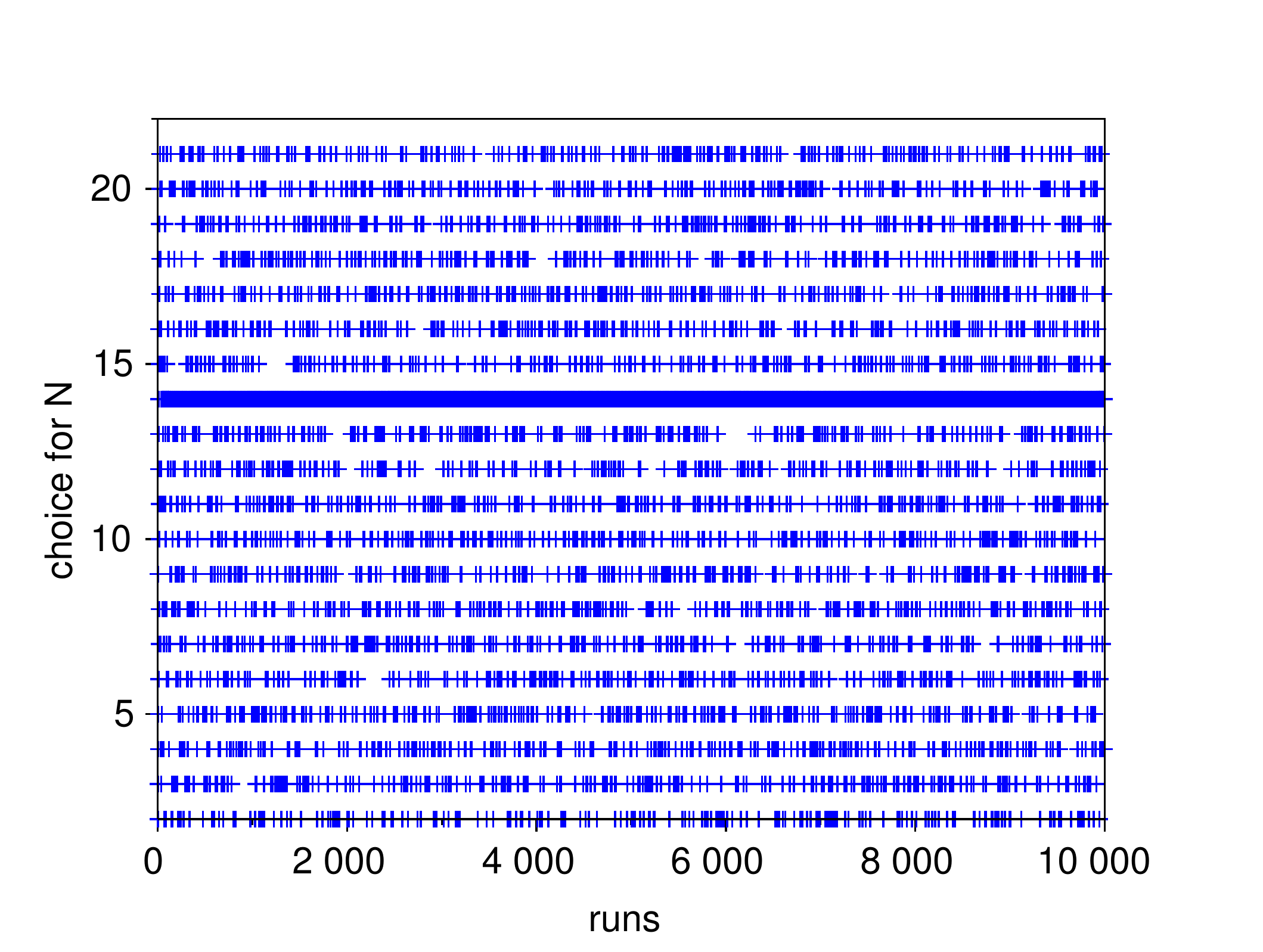}\includegraphics[width=7cm]{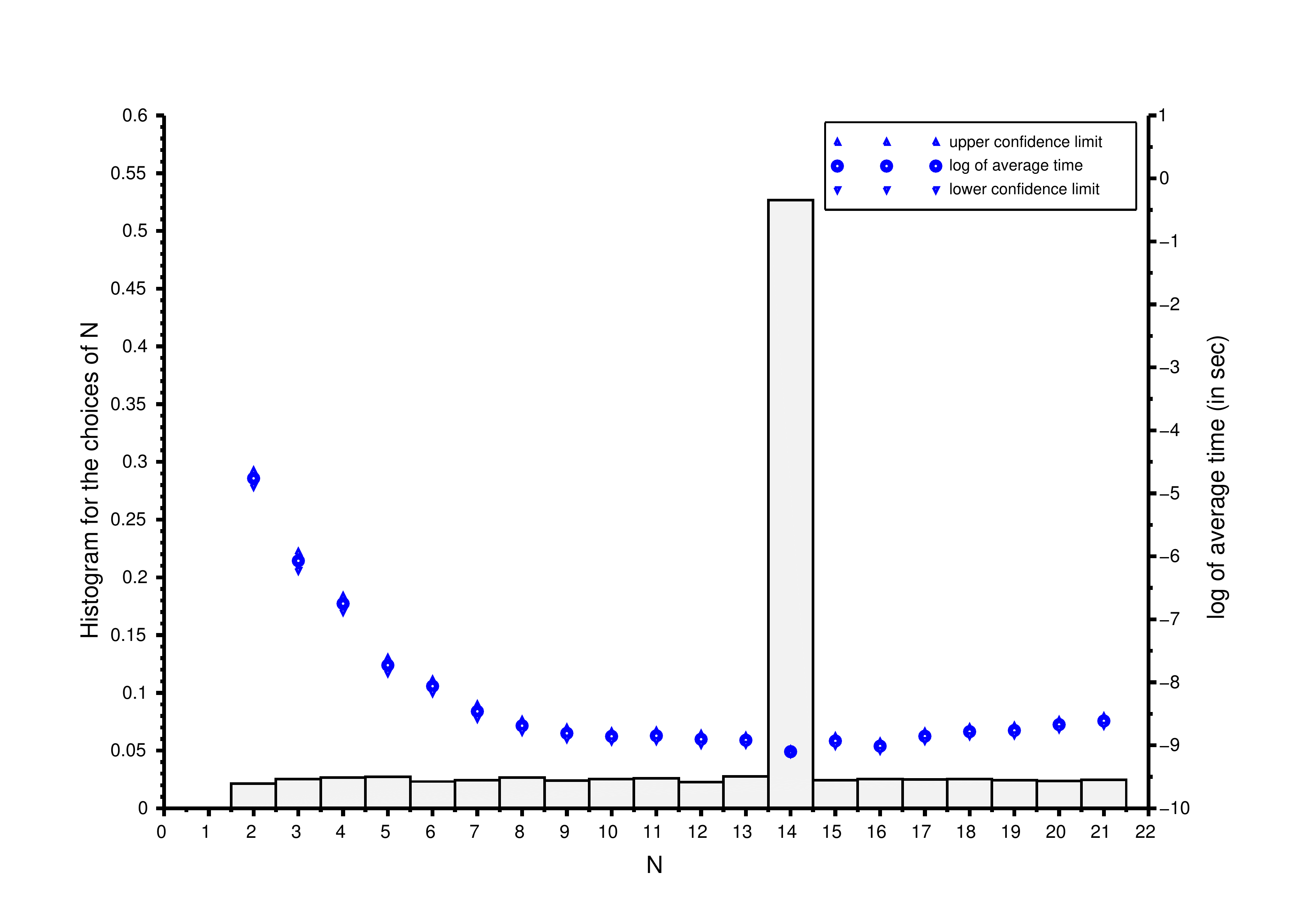}}
\centerline{\includegraphics[width=6cm]{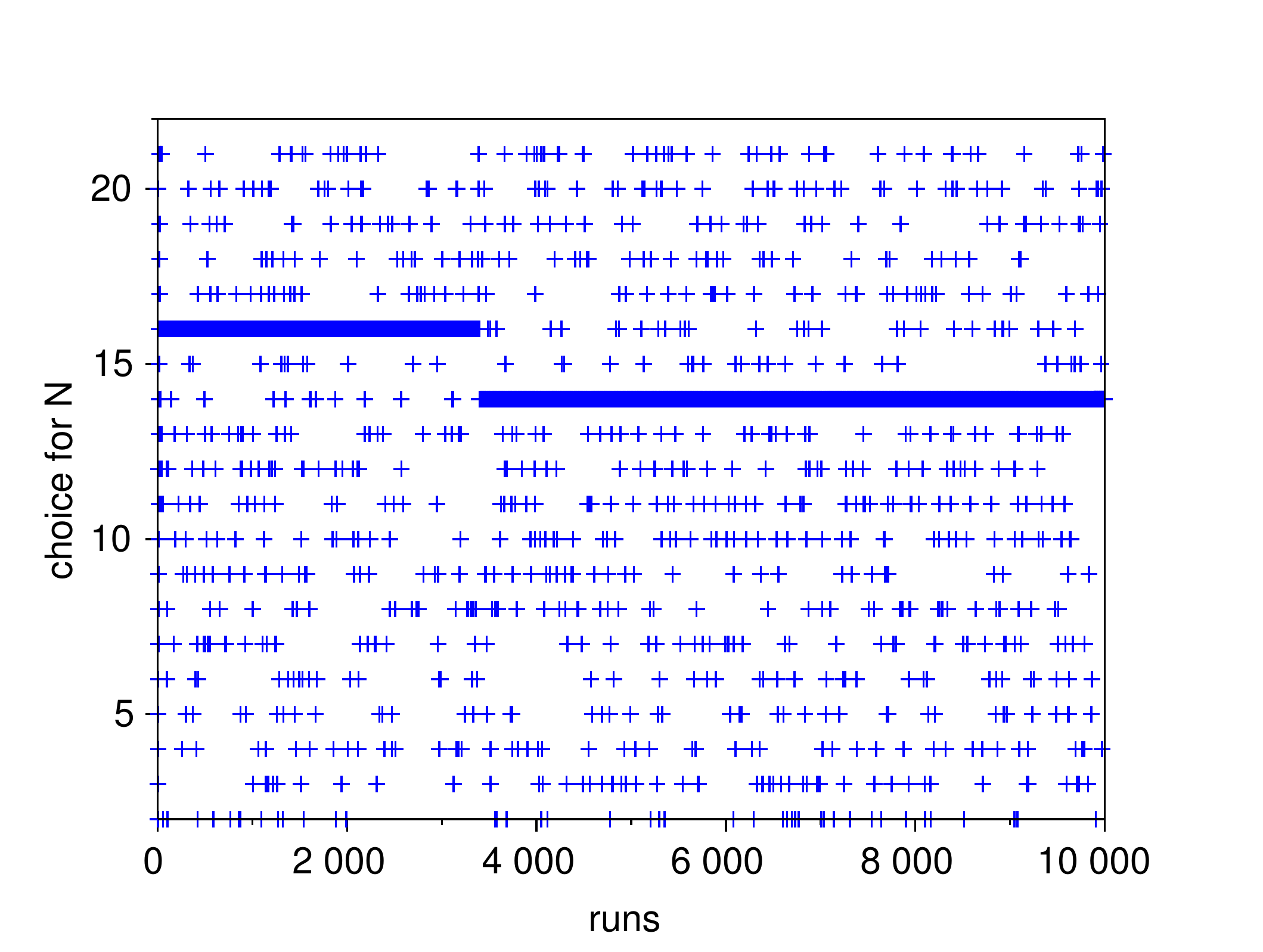}\includegraphics[width=7cm]{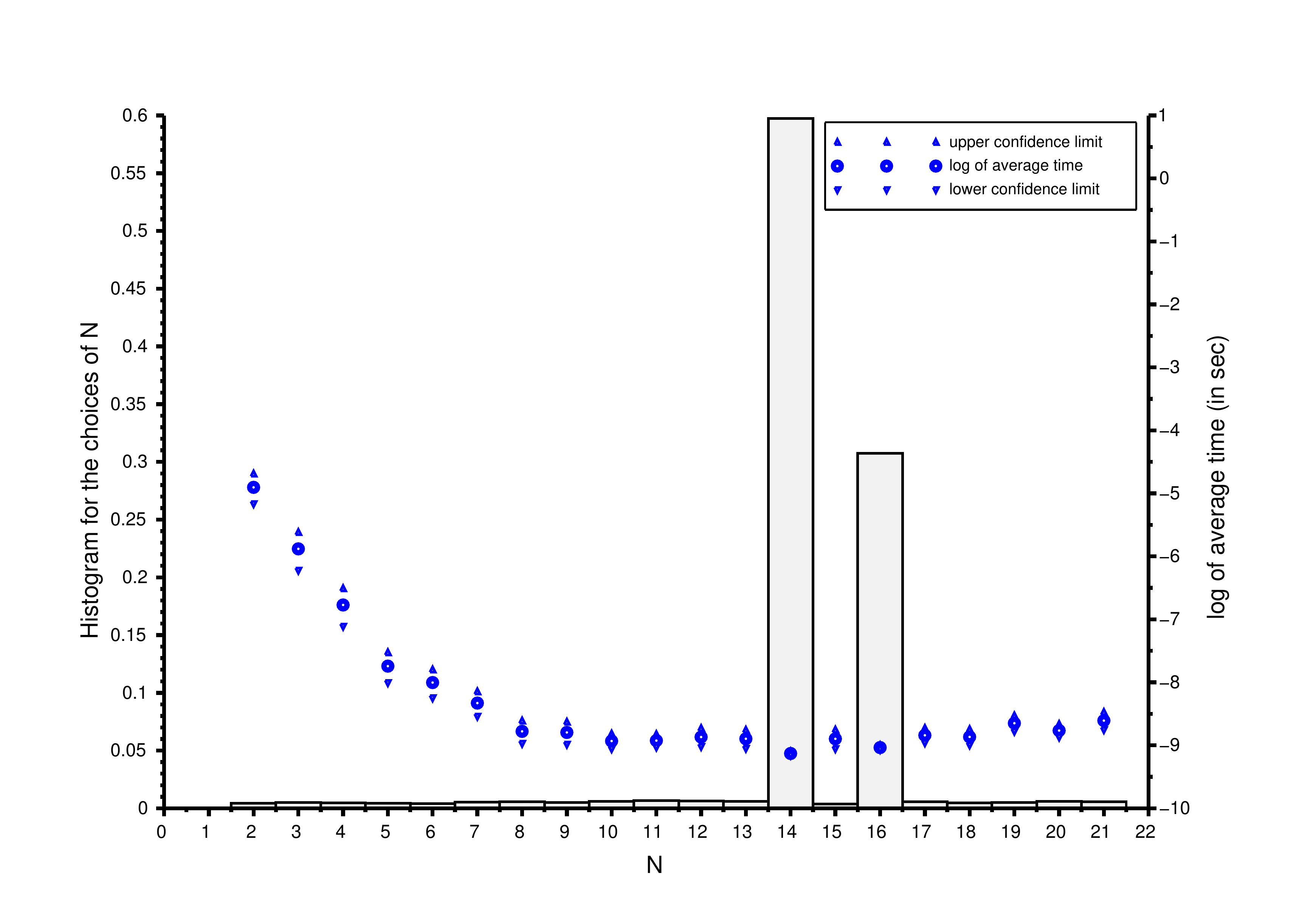}}
\caption{\small Choice of the box size (parameter $N$) versus the number of iterations (left) when the box size is chosen uniformly inbetween $2$ and $21$ accordingly to the $\epsilon$-greedy algorithm with $\epsilon=1$ (top), $\epsilon=0.5$ (middle) or $\epsilon=0.1$ (bottom) and histogram of the box size $N$ for a sample of size $10\,000$.  Here we consider the exit time of the interval $[a,b]=[0,7]$ and starting point $x=3$ and $T=1$. }
\label{Fig2n}
\end{figure}
\clearpage

The multi-armed bandit approach permits to possibly change the box size used for the exit time simulation by selecting the parameter $N$ inbetween a set of given values (here $\{2,\ldots, 21\}$, the arms of the bandit). The sequence of successive choices is randomized since the parameter $\epsilon$ which represents in some sense the level of noise (the proportion of exploration in the whole sequence of successive runs), belongs to $]0,1]$. In other words, the particular choice $\epsilon=1$ corresponds to a sequence of independent uniformly distributed choices whereas $\epsilon$ close to $0$ corresponds to a sequence of mainly deterministic choices linked to the argmin of the previous rewards (here the rewards are the consumption times).  

In Figure \ref{Fig2n}, we illustrate the behavior of the algorithm for three different values of $\epsilon$. In each case, the selections of the parameter $N$ throughout the sequence of iterations are represented by crosses in the figures (left). Once all exit times have been simulated, the assessment is represented by both the frequencies of each value of $N$ (histogram - right) and the corresponding average consumption time (with possibly its confidence interval). We can immediately observe the following.
\begin{itemize}
\item  In the case $\epsilon=1$, the choice of the parameter $N$ at each step of the algorithm does not depend on the previously observed consumption times, $N=14$ corresponding to its argmin is not privileged.
\item In the second case studied ($\epsilon=0.5$, middle), the particular choice  $N=14$ is rapidly privileged even if the relatively important level of noise implies a frequent visit of each proposed choice: $2,3,\ldots,21$. The exploration is quite important in that case.
\item Finally in the third case ($\epsilon=0.1$), the experiment leads to the following observation: the bandit algorithm makes $N=16$ its first choice but after a while (about 4\,000 iterations) the noise permits to leave this local minimum and to choose the global one.
\end{itemize}
So in order to reach a global minimum, it seems to be important not to choose the noise level $\epsilon$ too small. However we notice that the consumption times observed for both $N=14$ and $N=16$ are very close together, so the investigation of the argmin is not a crucial challenge.\medskip

\begin{figure}[h]
\centerline{\includegraphics[width=7cm]{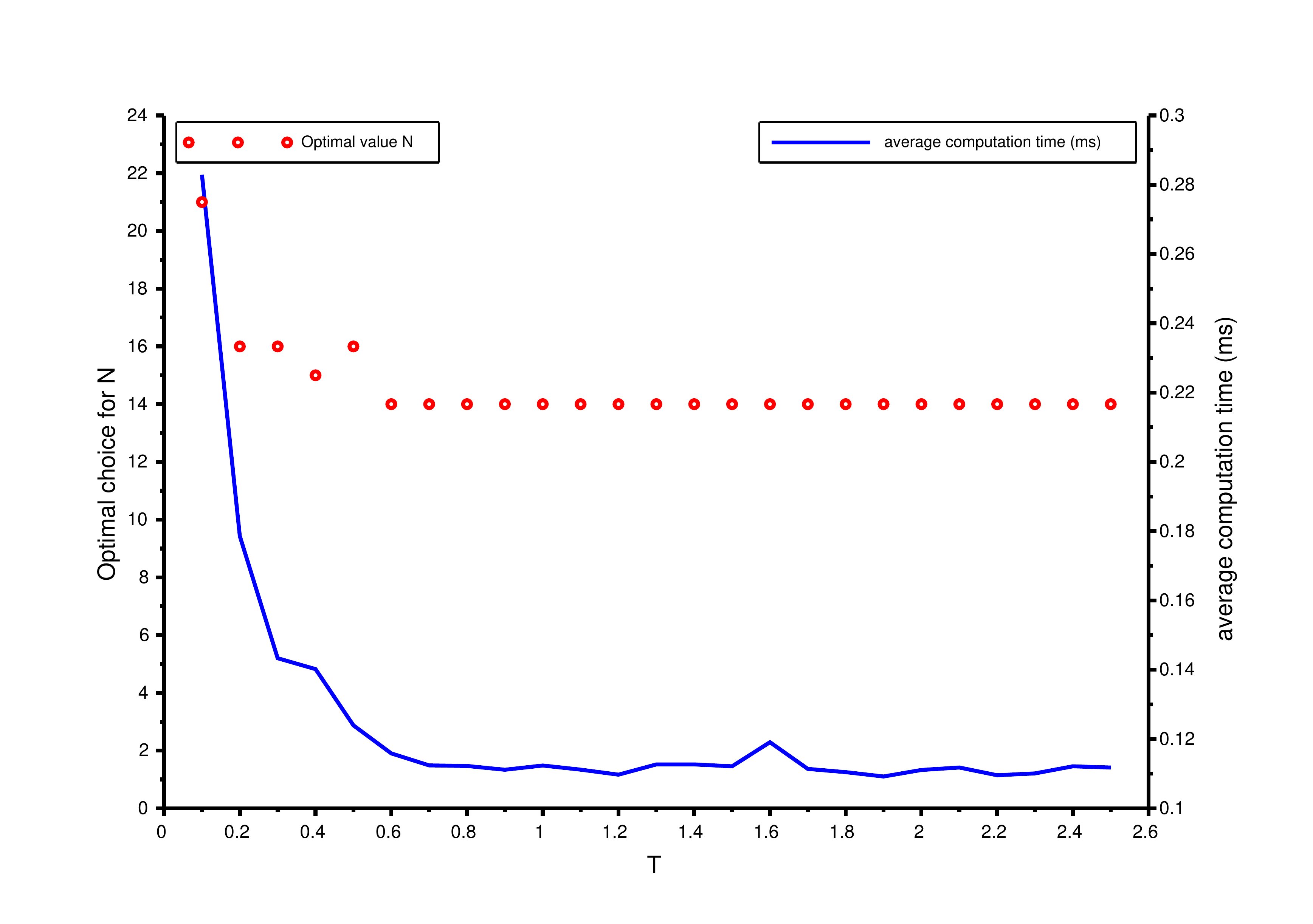}}
\caption{\small Optimal choice of the parameter $N$ and average time consumption versus $T$. We recall that the box size is $T\times (b-a)/N$. Here $N$ is chosen in the set $\{2,\ldots,21\}$ accordingly to the $\epsilon$-greedy algorithm with $\epsilon=0.1$, and the average is computed using a sample of size $10\,000$.  We consider the exit time from the interval $[a,b]=[0,7]$ and with starting value $x=3$. }
\label{fig3}
\end{figure}

Of course the box size of the basic components is the essential lever for the efficiency of the exit time simulation but it does not only depend on $N$. The area of the box is $2T\times (b-a)/N$ so that both $T$ and $N$ have to be correctly chosen. In Figure \ref{Fig2n}, the parameter $T$ is fixed ($T=1$) whereas $N$ varies. Once the optimal choice of $N$ is emphasized, it is possible to observe how it depends on $T$. Figure \ref{fig3} illustrates that the consumption time of the algorithm does actually not precisely depend on $T$, provided that $T$ is not too small ($T=1$ is a reasonable choice). 
%
%
%
%
%
%
%
%
%
%
%
%
%
%
%
%
%
%
%
%
%
%
%
%
\subsection{Ornstein-Uhlenbeck processes} 
Let us now consider a diffusion process which does not verify the particular condition presented in Remark \ref{rem:Tinfinite}. We aim to illustrate the efficiency of the bandit algorithm with the Ornstein-Uhlenbeck processes. So we consider the stochastic process with unitary diffusion coefficient and drift term  $\mu(x)=-\lambda x$. The aim is to simulate in some efficient way the first exit time of the interval  $[a,b]$. Since the process is mean reverting, its behavior will depend on the location of $0$, either in the interval $]a,b[$ either on the boundary or outside that interval. In order to present a complete illustration, we focus our attention on two different examples:
\begin{itemize}
\item Ex.1: interval $[a,b]=[0,7]$, drift $\lambda=1$ and starting position $x=3$.
\item Ex.2: interval $[a,b]=[-2,2]$, drift  $\lambda=2$ and starting position $x=0.5$.
\end{itemize}
In oder to simulate the first exit time from $[a,b]$ we aim to compare $\diff$ with the multi-armed bandit approach $\bdiff$. Let us just recall that $\diff(x,[a,b],T)$ is based on a sequential observation of the paths on the intervals $[nT,(n+1)T]$, $n\ge 0$, till the exit happens. Here $T$ is a parameter which influences the efficiency of the numerical procedure. For both cases under consideration, we observe that $T=0.5$ is a reasonable choice as suggested by the following table.  It presents the estimated computation times in \emph{ms} for one exit time generation (estimation with a sample of $10\,000$ exit times).
\\

%
%
%
%
%
%
%

%
%





\begin{figure}[h]
\centerline{\hspace*{0.5cm}\includegraphics[width=6.5cm]{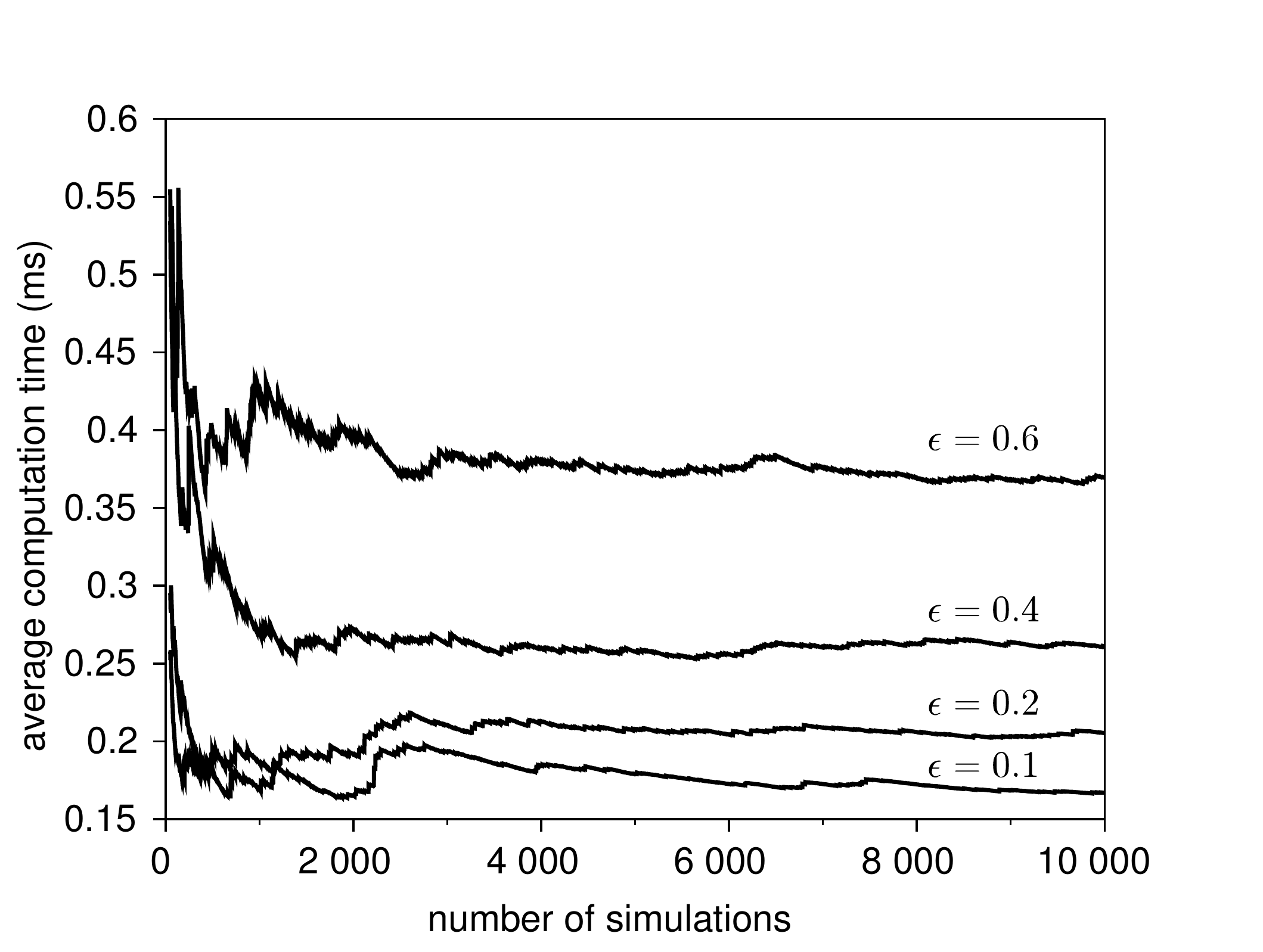}\hspace*{-0.5cm}\includegraphics[width=6.5cm]{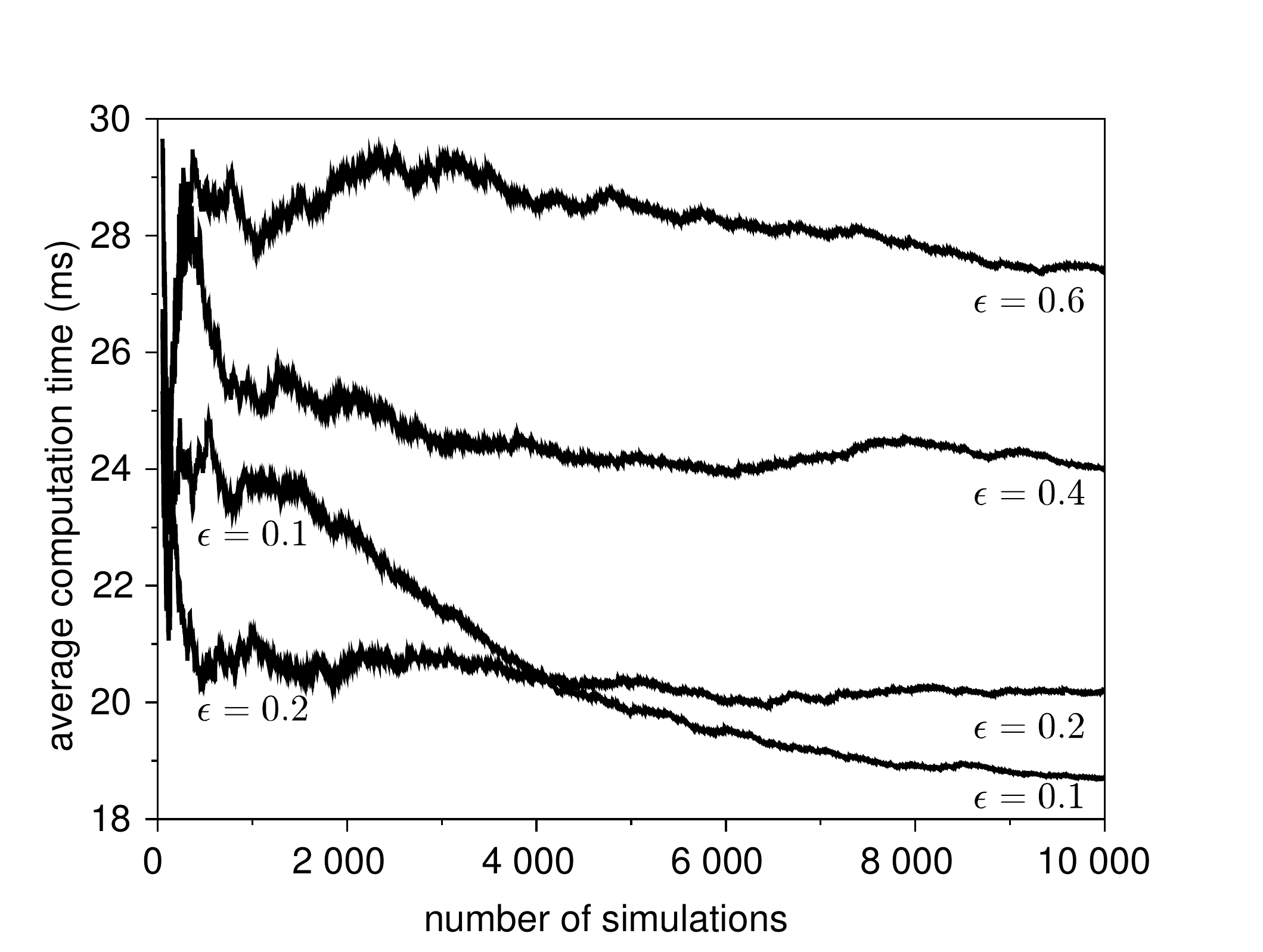}}
\caption{\small Empirical mean of the computation time (in \emph{ms}) versus the number of simulations of exit times from the interval $[a,b]$ for both examples (Ex.1 left and Ex.2 right).  We use different $\epsilon$-greedy algorithms with $T=0.5$ and $N$ is chosen in the set $\{ 2,\ldots, 21 \}$.}
\label{FigOUgreedy}
\end{figure}

\renewcommand{\arraystretch}{1.2}
\centerline{\begin{tabular}{|p{0.5cm}||p{1.2cm}|p{1.2cm}|p{1.2cm}||p{1.2cm}|p{1.2cm}|p{1.2cm}|}
\hline
 & \multicolumn{3}{|c||}{computation time in \emph{ms} (Ex.1)} & \multicolumn{3}{|c|}{computation time in \emph{ms}  (Ex.2)}\\ 
\hline
$T$ & average &  \multicolumn{2}{|c||}{confidence ($95\%$)
} & average &  \multicolumn{2}{|c|}{confidence ($95\%$)}\\
\hline
\hline
 $\ 0.1$ & $\ 4.720$ & $\ 4.602$ & $\ 4.837$ & $\ 68.636$ & $\  67.304 $ & $\ 69.969$  \\
$\ 0.2$ & $\ 4.292$ & $\ 4.194$ & $\ 4.389$ & $\ 52.355$ & $\ 51.346$ & $\ 53.363$ \\
\colorbox[gray]{0.85}{$0.5$} & \colorbox[gray]{0.85}{$4.195$} &  \colorbox[gray]{0.85}{$4.123$} &  \colorbox[gray]{0.85}{$4.267$} & \colorbox[gray]{0.85}{$46.471$} &  \colorbox[gray]{0.85}{$45.542$} &  \colorbox[gray]{0.85}{$47.401$} \\
$\ 1$ & $\ 5.022$ & $\ 4.929$ & $\ 5.115$ & $\ 55.076$ & $\ 53.973$ & $\ 56.178$\\
$\ 2$ & $\ 8.001$ & $\ 7.852$ & $\ 8.151$ & $\ 98.056$ & $\ 96.121$ & $\ 99.990$ \\
$\ 3$ & $13.135$ & $12.876$ & $ 13.394$ & $ 186.576$ & $ 182.906$ & $ 190.246$ \\
\hline
\end{tabular}}
\vspace*{0.2cm}
The consumption times associated with
 $T=0.5$ become therefore our reference values which need to be compared to the times issued from the multi-armed bandit. Figure \ref{FigOUgreedy} illustrates the efficiency of our approach for both examples (Ex.1 and Ex.2) since these consumption times have been reduced especially for small noise intensity $\epsilon$ (we suggest to choose $\epsilon$ smaller than $0.5$). This acceleration is less impressive when the origin $0$ belongs to the interval $[a,b]$ (Figure \ref{FigOUgreedy}, right). Let us also note that the curves of the average computation time associated with the parameters $\epsilon=0.1$ and $\epsilon=0.2$ intersect each other: if one needs a huge number of simulations, then one prefer $\epsilon=0.1$ which permits to find the global minimum and to avoid to often visit the other values of $N$. If one needs rather an intermediate value of simulations (for instance, $1\,000$ exit times), then it is better to increase a little the noise in the multi-armed bandit ($\epsilon=0.2$) in order to find quickly the optimal value $N$ even if the algorithm  frequently visits all the other values of $N$. 

\begin{figure}[h]
\centerline{\hspace*{0.5cm}\includegraphics[width=7cm]{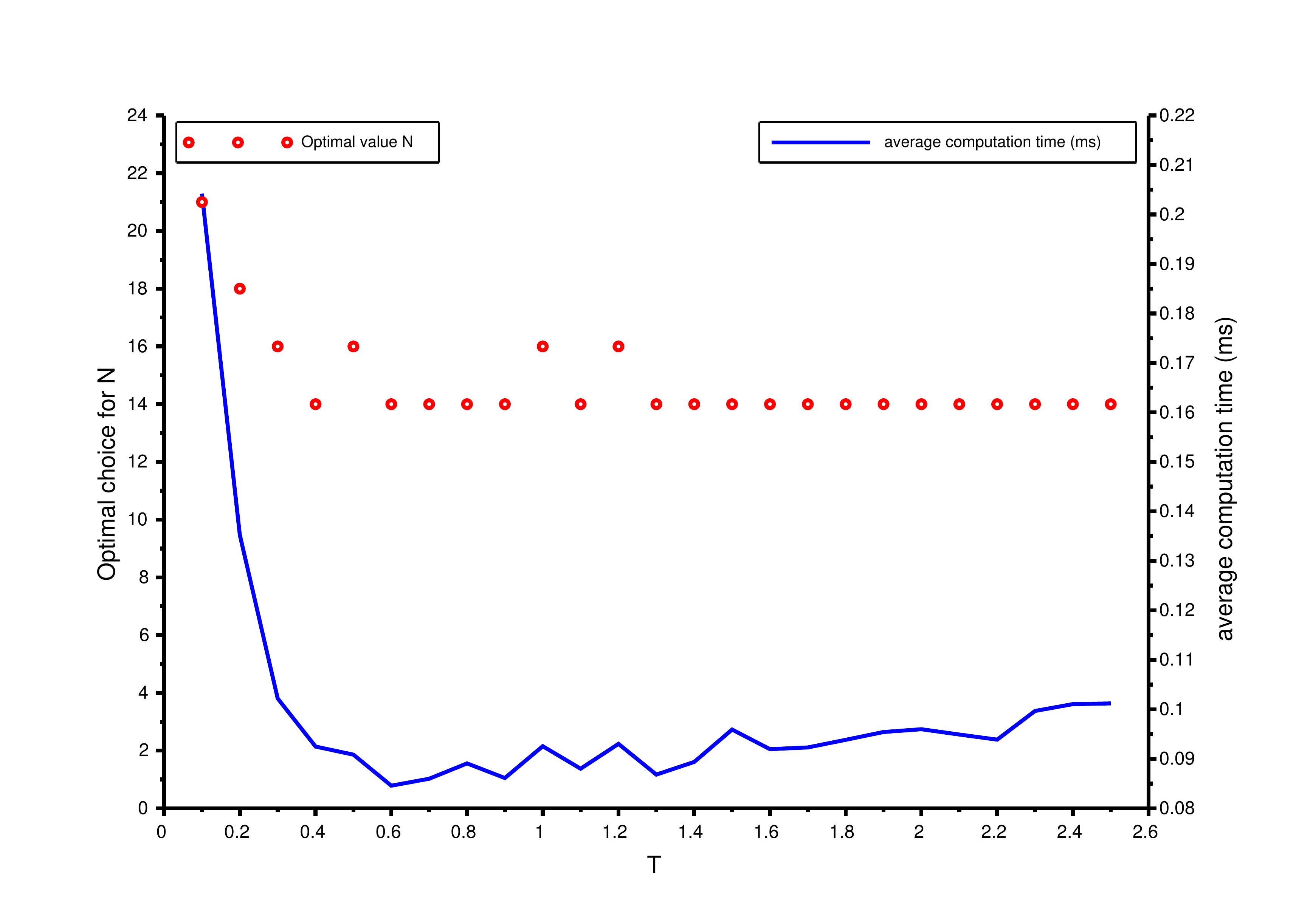}\hspace*{-0.5cm}\includegraphics[width=7cm]{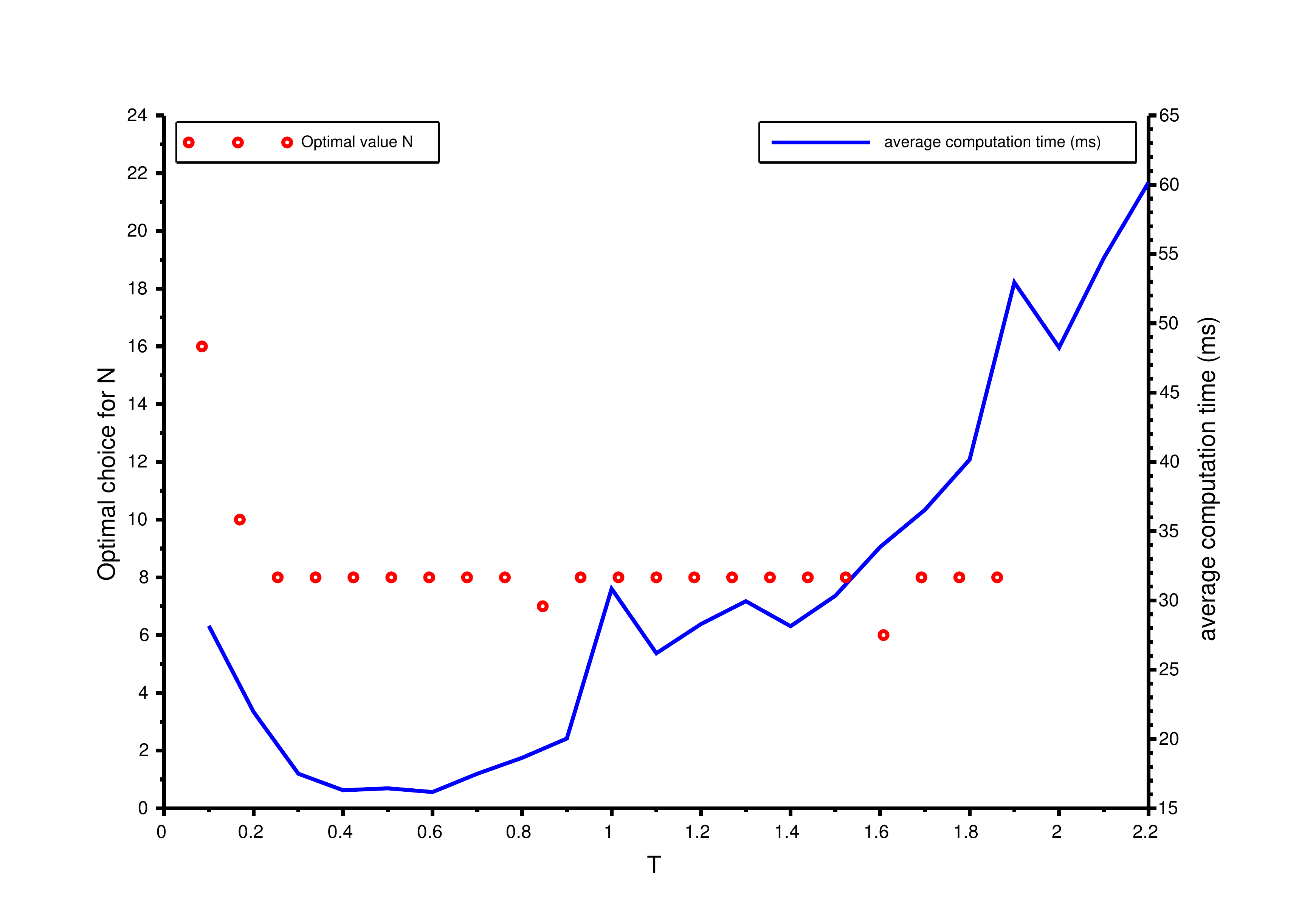}}
\caption{\small Optimal choice of the parameter $N$ and average time consumption versus $T$ for Ex.1 (left) and Ex.2 (right) . Here $N$ is chosen in the set $\{2,\ldots,21\}$ accordingly to the $\epsilon$-greedy algorithm with $\epsilon=0.1$, and the average is computed using a sample of size $10\,000$. }
\label{FigOUgreedy-T}
\end{figure}

The parameter $T$ was fixed so far in the study of the Ornstein-Uhlenbeck process and the attention was focused on the best choice of $N$. As already explained in Section \ref{sec:ex1}, the \bex \ algorithm depend both on $N$ and $T$.  Figure \ref{FigOUgreedy-T} represents the dependence of the optimal choice of $N$ and the average time consumption with respect to the parameter $T$. This illustration emphasizes that the efficiency is not strongly dependent with respect to $T$ provided that $T$ is neither too small nor too large. Even if the box size depends on both parameters $T$ and $N$,  it is therefore more clever to look after the best choice for $N$ rather than the best choice for $T$. Moreover we prefer to avoid an application of the $\epsilon$-greedy algorithm to the couple $(T,N)$ ($T$ would be discretized) trying to keep things simple.

\subsection{Cox-Ingersoll-Ross Processes}
In all the previous examples, the diffusions under observation have a unitary diffusion coefficient. In such situations, both \diff\ and \bdiff\ can be applied directly without using Lamperti's transform (see Remark \ref{rem:lamperti}). In order to complete the numerical illustration, we introduce a third example linked to the so-called CIR model (Cox-Ingersoll-Ross) which is of prime importance in the mathematical finance framework, in particular for the modelization of interest rates. The CIR model is characterized by the following stochastic differential equation:
\begin{equation}
\label{eq:CIR}
dX_t=k(\theta-X_t)\,dt+\sigma \sqrt{X_t}\,dB_t,\quad t\ge 0, \quad X_0=x>0.
\end{equation}
Here $k$ and $\theta$ are two parameters. Since the diffusion coefficient is not constant, we have to use the Lamperti transformation introduced in \eqref{eq:defofS}. So we define $\mathcal{S}(x)=\frac{2}{\sigma}\,\sqrt{x}$. Then $\widehat{X}_t:=\mathcal{S}(X_t)$ is a diffusion process with unitary diffusion coefficient and drift term:
\[
\mu_0(x)=\frac{\rho}{x}-\frac{kx}{2}\quad\mbox{where}\quad \rho:=\frac{(4k\theta-\sigma^2)}{2\sigma^2}.
\]
Let us assume that the parameters appearing in \eqref{eq:CIR} satisfy the condition: $\rho>0$. Consequently the CIR process starting from a positive initial point stays strictly positive (see, for instance, \cite{yor-jeanblanc} in Section 6.3.1) and the function $\gamma$ and $\beta$ used in the algorithms have an explicit expression easy to handle with:
\begin{equation}
\gamma(x)=\frac{1}{2}\left( \left(\frac{\rho}{x} -\frac{kx}{2}\right)^2 -\frac{\rho}{x^2}-\frac{k}{2} \right),\quad 
\beta(x)=x^\rho e^{-k x^2/4}.
\end{equation}
For numerical illustration, we deal with the exit problem from $[a,b]=[1,6]$ for the CIR model starting in $x=3$ with coefficients $k=3$, $\theta=7$ and $\sigma=1$. As in the Ornstein-Uhlenbeck context, we can here use the \bex (x,[$\mathcal{S}(a),\mathcal{S}(b)]$,T)\ algorithm in order to simulate both the exit time and the exit location. This algorithm depends on a parameter $T$. We obtain the following average computation times for one exit time generation:\\

\renewcommand{\arraystretch}{1.2}
\centerline{\begin{tabular}{|p{0.5cm}||p{2.5cm}|p{1.7cm}|p{1.7cm}|}
\hline
$T$ & average in \emph{ms} &  \multicolumn{2}{|c|}{confidence ($95\%$)
} \\
\hline
\hline
 $\ 0.1$ & $\ 12.206$ & $\ 11.981$ & $\ 12.431$  \\
$\ 0.2$ & $\ 11.939$ & $\ 11.718$ & $\ 12.161$  \\
$\ 0.5$ & $\ 11.901$ &  $\ 11.671$ &  $\ 12.130$ \\
$\ 1$ & $\ 12.272$ & $\ 12.036$ & $\ 12.509$ \\
\hline
\end{tabular}}

\begin{figure}[h]
\centerline{\hspace*{0.5cm}\includegraphics[width=6.5cm]{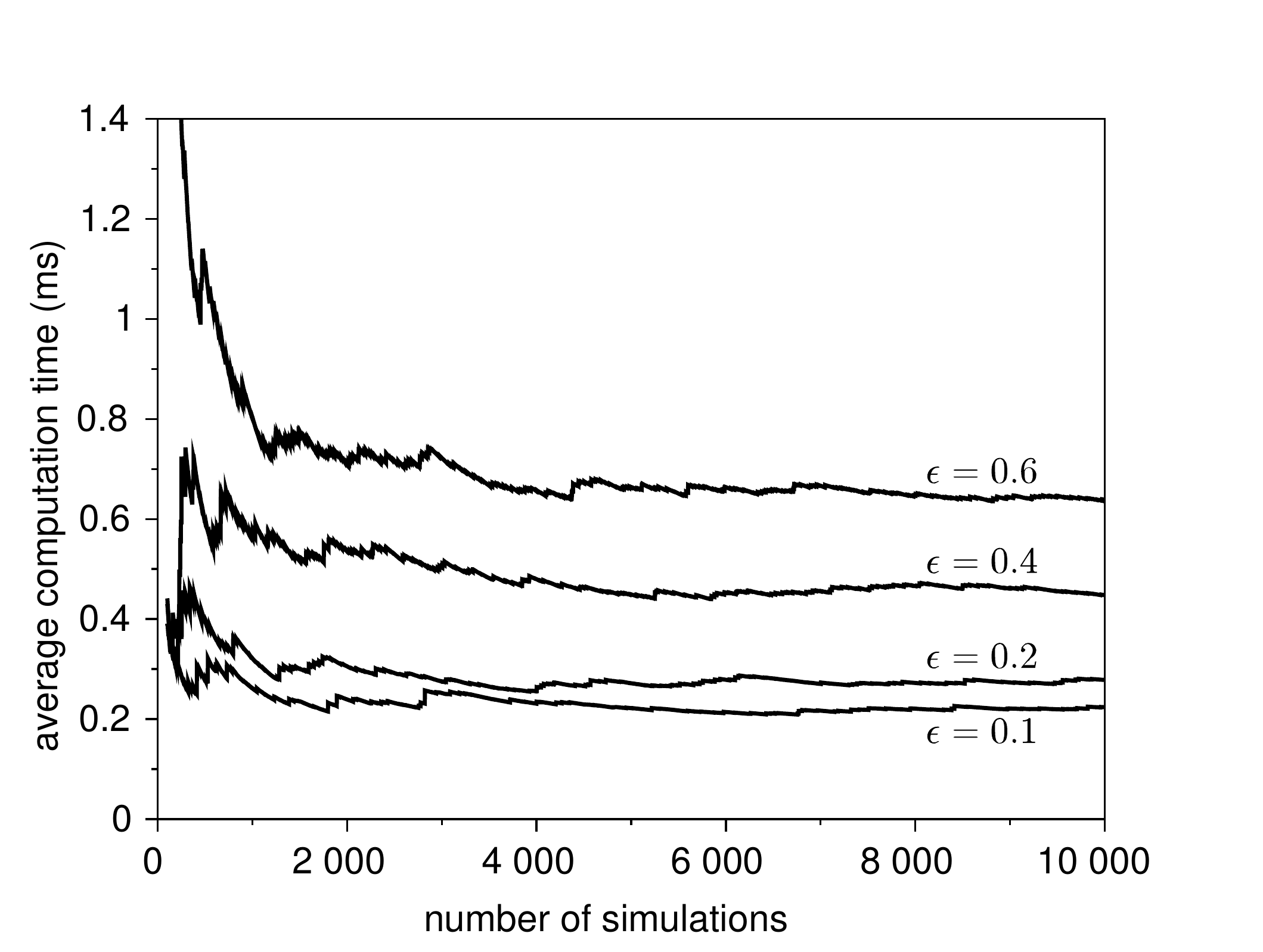}\hspace*{-0.5cm}\includegraphics[width=6.5cm]{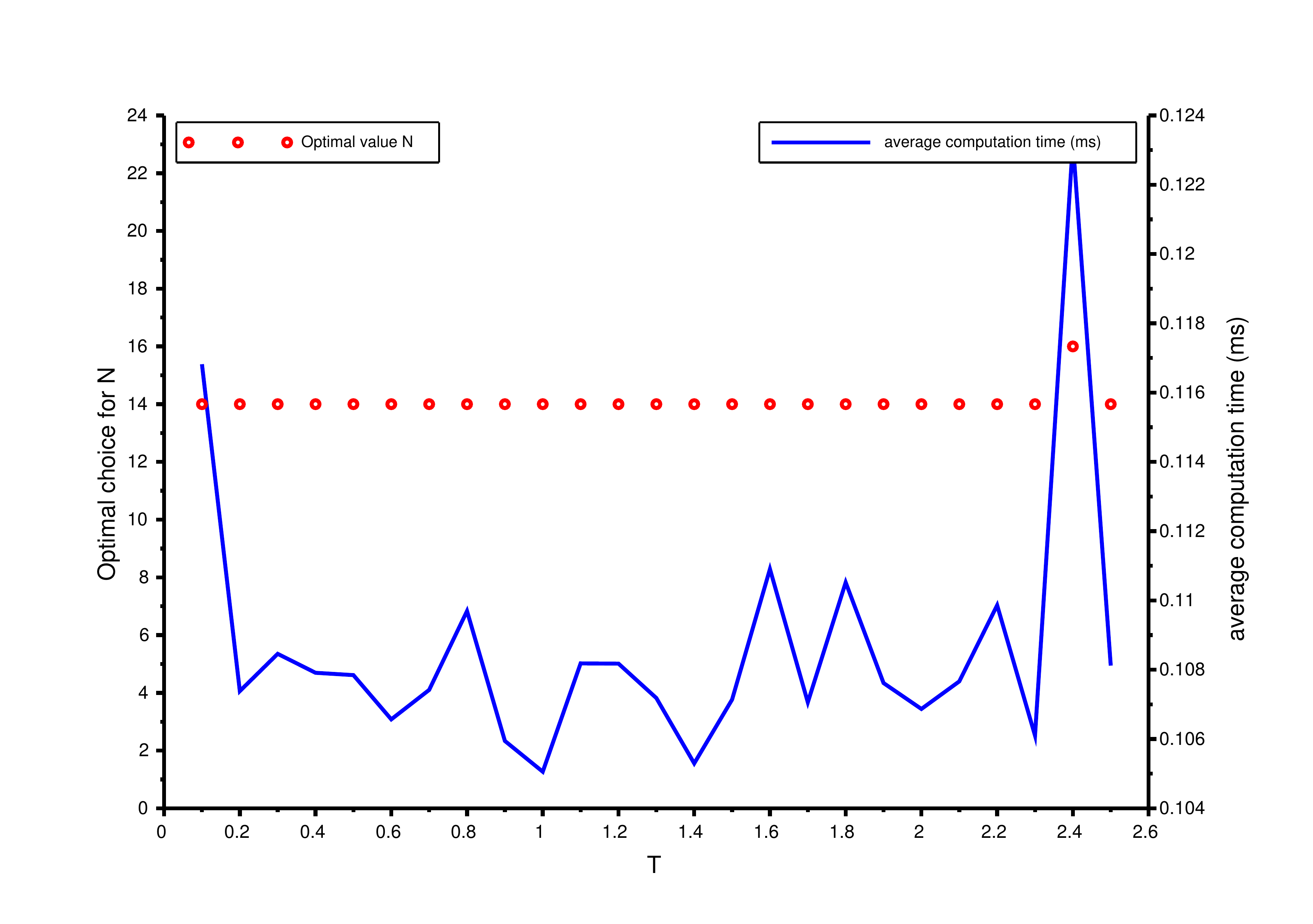}}
\caption{\small Empirical mean of the computation time (in \emph{ms}) versus the number of simulations of exit times from the interval $[1,6]$ with different $\epsilon$-greedy algorithms and $T=0.5$ (left). Optimal choice of the parameter $N$ and average time consumption versus $T$ for the CIR model (right). Here $N$ is chosen in the set $\{2,\ldots,21\}$ accordingly to the $\epsilon$-greedy algorithm with $\epsilon=0.1$, and the average is computed using a sample of size $10\,000$.}
\label{FigCIRgreedy}
\end{figure}

\vspace*{0.3cm}
We can observe that the parameter $T$ has only a weak influence on the $\bex$ efficiency provided $T$ belongs to an interval of reasonable values (here inbetween $0.1$ and $1$). Let us now compare these computation times of the order of 12 ms per simulation to the \bdiff\ algorithm one. In Figure \ref{FigCIRgreedy} (left), we observe a significant time reduction as soon as $\epsilon$ (the parameter of the $\epsilon$-greedy procedure) is sufficiently small, we reach a computation time near to $0.2$ ms per simulation (for a sample size $10\,000$). Since the box size used in \bdiff\ depends on both parameters $N$ and $T$, we wonder if the optimal value of $N$ strongly depends on $T$. As we can see in Figure \ref{FigCIRgreedy} (right), it is not the case: there is neither large swings in the optimal choice of the value of $N$ nor in the average consumption time associated with this optimal $N$.

\subsection*{Conclusion}
The exact simulation procedure $\bex(x,[a,b],T)$ proposed in \cite{Herrmann-Zucca-2} permits to generate the exit time and exit location from an interval $[a,b]$ in the diffusion context. In this study, we emphasize a reinforcement learning method based on a multi-armed bandit which permits to accelerate the $\bex$ algorithm in any case. As presented in Section \ref{sec:illu}, sometimes the consumption time reduction is very strong and sometimes sensible. The tremendous advantage of the algorithm $\bdiff$ is its universality: it does not depend on the the particular family of diffusion under consideration.\medskip

Let us also note that the authors have chosen the $\epsilon$-greedy algorithm for the acceleration procedure since it is simple to explain and particularly efficient. Of course any other algorithm used in the classical multi-armed bandit problem can be tested for the acceleration of $\bex$.\medskip

\begin{framed}{\sl \noindent All the numerical tests have been done on the same computer:\\ Intel Core i5, 1.6 GHz}
\end{framed}
%
%
%
%

%
%
\end{document}